\documentclass[MSNbibl,number,citesort,dvips]{arxbj}
\usepackage{dcolumn,upgreek}
\usepackage[left]{vmkcol}
\usepackage{graphicx}

\aid{0}
\volume{19}
\issue{5B}
\pubyear{2013}
\firstpage{2557}
\lastpage{2589}
\doi{10.3150/12-BEJ466} %kopijuoti is PTS

\makeatletter

\newcommand{\rrvert}{\vert}
\newcommand{\llvert}{\vert}
\newcolumntype{d}[1]{D{.}{.}{#1}}
\newcommand{\eqref}[1]{(\ref{#1})}

\newcommand{\baF}{\bar{F}}
\newcommand{\habaF}{\hat{\hspace*{-2pt}\bar{F}}_n}
\newcommand{\toP}{\stackrel{P}{\longrightarrow}}
\newcommand{\tod}{\stackrel{d}{\longrightarrow}}
\newcommand{\R}{\mathbb{R}}
\newcommand{\PP}{\mathbb{P}}
\newcommand{\I}{\mathbb{I}}
\newcommand{\indic}{\mathbb{I}}
\newcommand{\E}{\mathbb{E}}

\newtheorem{Theo}{Theorem}
\newtheorem{prop}{Proposition}
\newtheorem{Coro}{Corollary}
\newtheorem{Lem}{Lemma}
\newremark{Remark}{Remark}

\makeatother

\begin{document}
\begin{frontmatter}

\title{On kernel smoothing for extremal quantile regression}
\runtitle{Extremal quantile regression}

\begin{aug}
%%%% inicialai - be tarpu
\author[1]{\fnms{Abdelaati} \snm{Daouia}\thanksref{1}},
\author[2]{\fnms{Laurent} \snm{Gardes}\thanksref{2}} \and
\author[3]{\fnms{St\'ephane} \snm{Girard}\thanksref{3}\corref{}\ead[label=e3]{Stephane.Girard@inria.fr}}
\runauthor{A. Daouia, L. Gardes and S. Girard} %% auto
\address[1]{Institute of Statistics, Catholic University of Louvain,
Belgium.
Toulouse School of Economics, University of Toulouse, France} %
\address[2]{Universit\'e de Strasbourg \& CNRS, IRMA, UMR 7501, France}
%%\printead{e2}
\address[3]{INRIA Rh\^one-Alpes \& LJK, team Mistis, Inovall\'ee, 655,
av. de l'Europe, Montbonnot, 38334 Saint-Ismier cedex, France.
\printead{e3}}
\end{aug}

% HISTORY:
\received{\smonth{10} \syear{2011}}
\revised{\smonth{7} \syear{2012}}

% ABSTRACT
%
\begin{abstract}
Nonparametric regression quantiles obtained by inverting a kernel
estimator of
the conditional distribution of the response are long established in
statistics.
Attention has been, however, restricted to ordinary quantiles staying away
from the tails of the conditional distribution.
The purpose of this paper is to extend their asymptotic theory far
enough into the tails. We focus on extremal quantile regression
estimators of
a response variable given a vector of covariates in the general setting,
whether the conditional extreme-value index is positive, negative, or zero.
Specifically, we elucidate their limit distributions when they are located
in the range of the data or near and even beyond the sample boundary, under
technical conditions that link the speed of convergence of their (intermediate
or extreme) order with the oscillations of the quantile function and a
von-Mises property of the conditional distribution.
A simulation experiment and an illustration on real data were presented.
The real data are the American electric data where the estimation of
conditional extremes is found to be of genuine interest.
\end{abstract}

% KEYWORDS
%

\end{frontmatter}

%s1 #&#
\section{Introduction}
Quantile regression plays a fundamental role in various statistical
applications. It complements the classical regression on the
conditional mean
by offering a more useful tool for examining how a vector of regressors
$X\in\mathbb{R}^p$ influences the entire distribution of a response variable
$Y\in\mathbb{R}$. The nonparametric regression quantiles obtained by
inverting a kernel
estimator of the conditional distribution function are used widely in applied
work and investigated extensively in theoretical statistics. See, for
example~\cite{BER,SAM,STO,STU}, among others. Attention has been,
however, restricted to
conditional quantiles having a fixed order $\alpha\in(0,1)$.
In the following, the order $\alpha$ has to be understood as the
conditional probability to be
larger than the conditional quantile. In result, the available large
sample theory does not apply sufficiently far in the tails.

There are many important applications in ecology, climatology,
demography, biostatistics, econometrics, finance, insurance, to name a few,
where extending
that conventional asymptotic theory further into the tails of the conditional
distribution is an especially welcome development.
This translates into considering the order $\alpha=\alpha_n\to0$ or
$\alpha_n\to1$ as the sample size $n$ goes to infinity.
Motivating examples include the study of extreme rainfall as a function
of the
geographical location~\cite{Gar}, the estimation of factors of
high risk in finance~\cite{Tsy}, the assessment of the optimal cost of the
delivery activity of postal services~\cite{DFS},
the analysis of survival at extreme durations~\cite{KON},
the edge estimation in image reconstruction~\cite{PAR},
the accurate description of the upper tail of the claim size
distribution for
reinsurers~\cite{BES},
the analysis of environmental time series with application to trend detection
in ground-level ozone~\cite{SMI}, the estimation of autoregressive models
with asymmetric innovations~\cite{FEIRES}, etc.

There have been several efforts to treat the asymptotics of extreme
conditional quantile estimators in semi/parametric and other
nonparametric regression models.
For example, Chernozhukov~\cite{CHE} and Jureckov\'a~\cite{JUR}
considered the extreme
quantiles in the linear regression model and derived their asymptotic
distributions under various distributions of errors.
Other parametric models are proposed in~\cite{DAV,SMI}, where some
extreme-value based techniques are extended to the
point-process view of high-level exceedances.
A~semi-parametric approach to modeling trends in sample extremes, based on
local polynomial fitting of the Generalized extreme-value distribution,
has been introduced in~\cite{DAR}.
Hall and Tajvidi~\cite{HAL} suggested a nonparametric
estimation of the temporal trend when fitting parametric models to extreme
values. Another semi-parametric method has been developed in~\cite{BEG}, where the regression is based on a Pareto-type conditional
distribution of the response.
Fully nonparametric estimators of extreme conditional quantiles have been
discussed in~\cite{BEG,CHA}, where the former approach is based on the
technique of local
polynomial maximum likelihood estimation, while spline estimators are fitted
in the latter by a maximum penalized likelihood method.
Recently, \cite{Gar2,GGL} proposed,
respectively, a moving-window based
estimator for the tail index and extreme quantiles of heavy-tailed
conditional distributions, and
they established their asymptotic properties.

In the context of kernel-smoothing, the asymptotic theory
for quantile regression in the tails is relatively unexplored and still in
full development. Daouia \textit{et al.}~\cite{Test} have extended the
asymptotics further into the tails in the particular
setting of a heavy-tailed conditional distribution, while~\cite{Gir,Gir2} have analyzed the case $\alpha_n=1/n$ in the particular
situation where the response $Y$ given $X=x$ is uniformly distributed.
The purpose of this paper is to develop a unified asymptotic
theory for the kernel-smoothed conditional extremes in the general setting
where the conditional distribution can be short, light or heavy-tailed.
We will focus on the $\alpha_n\to0$ case, which corresponds to the
class of
large quantiles of the upper conditional tail. Similar considerations
evidently apply to the case $\alpha_n\to1$. Specifically, we first obtain
the asymptotic normality of the extremal quantile regression
under the `intermediate' order condition $nh^p\alpha_n\to\infty$
where $h=h_n\to0$ stands for the bandwidth involved
in the kernel smoothing estimation.
Next, we extend the asymptotic normality far enough into the `most extreme'
order-$\beta_n$ regression quantiles with $\beta_n/\alpha_n\to0$, thus
providing a conditional analog of modern extreme-value results~\cite{deHaanFer}. We also analyze kernel-smoothed Pickands type estimators
of the conditional extreme-value index as in the familiar nonregression
case~\cite{Drees}.

The paper is organized as follows. Section~\ref{notation} contains the basic
notation and assumptions. Section~\ref{results} states the main results.
Section~\ref{mc} presents some simulation evidence and practical guidelines.
Section~\ref{data} provides a motivating example in production theory, and
the \hyperref[appendix]{Appendix} collects the proofs.

%s2 #&#
\section{The setting and assumptions}\label{notation}
%%%%%%%%%%%%%%%%%%%%%%%%%%%%%%%%%%%%%%%%%%%%%%%%%%%%%%%%%%%%%%%%%%%%%%%%%%%%%%%%%%%%%%%%%%%%%%%%%

Let $(X_i,Y_i)$, $i=1,\ldots,n$, be independent copies of a random pair
$(X,Y) \in\R^p \times\R$. The conditional survival function (c.s.f.) of
$Y$ given $X=x$ is
denoted by $\baF(y|x)=\PP(Y>y|X=x)$ and the probability density
function (p.d.f.) of $X$ is denoted by $g$.
We address the problem of estimating extreme conditional quantiles
\[
q(\alpha_n|x) = \baF^{\leftarrow}(\alpha_n|x)=\inf
\bigl\{ t,  \baF(t|x) \leq\alpha_n \bigr\},
\]
where $\alpha_n \to0$ as $n$ goes to infinity. In the following, we
denote by $y_F(x) = q(0|x) \in(-\infty,\infty]$ the endpoint of the
conditional distribution of $Y$ given $X=x$. The kernel estimator of
$\baF(y|x)$ is defined for all $(x,y)\in\R^p \times\R$
by
%
%e1 #&#
\begin{equation}
\label{defestproba} \habaF(y|x)= \sum_{i=1}^n
K_h(x-X_i) \I\{Y_i>y\} \Big/ \sum
_{i=1}^n K_h(x-X_i) ,
\end{equation}
where $\I\{\cdot\}$ is the indicator function and $h=h_n$ is a nonrandom
sequence such that
$h\to0$ as $n\to\infty$.
We have also introduced $K_h(t)=K(t/h)/h^p$
where $K$ is a p.d.f. on $\R^p$.
In this context, $h$ is called the window-width.
Similarly, the kernel estimators of conditional quantiles
$q(\alpha|x)$
are defined via the generalized inverse of $\habaF(\cdot|x)$:
%
%e2 #&#
\begin{equation}
\label{defestquant} \hat q_n(\alpha| x ) = \hat{\hspace*{-2pt}\bar
{F}} {}^{\leftarrow}_n(\alpha|x)=\inf \bigl\{t, \habaF(t|x)\leq
\alpha \bigr\}
\end{equation}
for all $\alpha\in(0,1)$.
Many papers are dedicated to the asymptotic properties of this type of
estimator for fixed $\alpha\in(0,1)$: weak and strong consistency are
proved, respectively, in~\cite{STO} and~\cite{GAN}, asymptotic normality
being established in~\cite{STU,SAM,BER}.
In Theorem~\ref{thquant} below,
the asymptotic distribution of~(\ref{defestquant})
is investigated when estimating extreme quantiles, {that is},
when $\alpha=\alpha_n$ goes to 0 as the sample size $n$ goes to infinity.
The asymptotic behavior of such estimators then depends on the nature of
the conditional distribution tail. In this paper, we assume that the
c.s.f. satisfies the following von-Mises condition, see for instance~\cite{deHaanFer}, equation (1.11.30):

\begin{longlist}[(A.1)]
\item[(A.1)] The function $\baF(\cdot|x)$ is twice differentiable and
\[
\lim_{y \uparrow y_F(x)} \frac{\baF(y|x) \baF''(y|x)}{(\baF')^2(y|x)} = \gamma(x)+1,
\]
where $\baF'(\cdot|x)$ and $\baF''(\cdot|x)$ are, respectively, the
first and the
second derivatives of $\baF(\cdot|x)$.
\end{longlist}
Here, $\gamma(\cdot)$ is an unknown function of the covariate~$x$ referred
to as the conditional extreme-value index.
Let us consider, for all $z\in\R$, the classical $K_z$ function defined
for all $u\in\R$ by
\[
K_z(u)=\int_1^u v^{z-1}
\,\mathrm{d}v.
\]
The associated inverse function is denoted by $K_z^{-1}$. Then, (A.1)
~implies that there exists a positive auxiliary function $a(\cdot|x)$
such that,
%
%e3 #&#
\begin{equation}
\label{Fu} \lim_{y \uparrow y_F(x)} \frac{\baF
(y+t(x)a(y|x)|x)}{\baF
(y|x)}=\frac
{1}{K^{-1}_{\gamma(x)}(t(x))},
\end{equation}
where $t(x) \in\R$ is such that $1+t(x)\gamma(x) >0$.
Besides, (\ref{Fu})~implies in turn that the conditional distribution
of $Y$ given $X=x$ is in the maximum domain of attraction (MDA) of the
extreme-value distribution
with shape parameter $\gamma(x)$, see~\cite{deHaanFer}, Theorem~1.1.8,
for a proof.
The case $\gamma(x)>0$ corresponds to the Fr\'echet MDA and $\baF
(\cdot|x)$
is heavy-tailed
while
the case $\gamma(x)=0$ corresponds to the Gumbel MDA and $\baF(\cdot
|x)$ is
light-tailed.
The case $\gamma(x)<0$ represents most of the situations where $\baF
(\cdot|x)$ is short-tailed,
that is, $\baF(\cdot|x)$ has a finite endpoint $y_F(x)$, this is
referred to as the Weibull MDA.

The convergence~(\ref{Fu}) is also equivalent to
%
%e4 #&#
\begin{equation}
\label{F1quant} b(t,\alpha|x):=\frac{q(t\alpha|x)-q(\alpha
|x)}{a(q(\alpha|x)|x)} - K_{\gamma(x)}(1/t) \to0
\end{equation}
for all $t>0$ as $\alpha\to0$, see~\cite{deHaanFer}, Theorem~1.1.6.
For all $(x,x')\in\R^p\times\R^p$, the Euclidean distance between
$x$ and $x'$ is denoted by $d(x,x')$. The following Lipschitz condition
is introduced:

\begin{longlist}[(A.2)]
\item[(A.2)] There exists $c_g>0$ such that
$ \llvert g(x)-g(x')\rrvert
\leq c_g d(x,x')$.
\end{longlist}

The last assumption is standard in the kernel estimation framework.

\begin{longlist}[(A.3)]
\item[(A.3)] $K$ is a bounded p.d.f. on $\R^p$, with support $S$
included in
the unit ball of $\R^p$.
\end{longlist}

%s3 #&#
\section{Main results}\label{results}
%%%%%%%%%%%%%%%%%%%%%%%%%%%%%%%%%%%%%%%%%%%%%%%%%%%%%%%%%%%%%%%%%%%%%%%%%%%%%%%%%%%%%%%%%%%%%%%%%%%%%%%%%%%%%%

Let $B(x,h)$ be the ball centered at $x$ with radius $h$. The
oscillations of the c.s.f. are controlled by
\[
\Delta_\kappa(x,\alpha) := \sup_{(x',\beta)\in B(x,h)\times
[\kappa
\alpha,\alpha]} \biggl\llvert
\frac{\baF(q(\beta|x)|x')}{\beta}-1 \biggr\rrvert,
\]
where $(\kappa,\alpha)\in(0,1)^2$. Under assumption (A.1), $\baF
(\cdot|x)$ is
differentiable
and the associated conditional density will be denoted in the sequel
by $f(\cdot|x)$. We first establish the asymptotic normality of $\hat
q_n(\alpha_n|x)$.

%th1 #&#
\begin{Theo}
\label{thquant}
Suppose \textup{(A.1)}, \textup{(A.2)} and \textup{(A.3)} hold. Let $0<\tau_J< \cdots< \tau_2 <
\tau_1
\leq1$ where $J$ is a positive integer
and $x\in\R^p$ such that $g(x)>0$.
If $\alpha_n \to0$ and there exists $\kappa\in(0,\tau_J)$ such that
\[
nh^p\alpha_n \to\infty, \qquad nh^{p}
\alpha_n \bigl(h \vee\Delta_{\kappa
}(x,\alpha_n)
\bigr)^2 \to0,
\]
then, the random vector
\[
\Bigl\{ f \bigl(q(\alpha_n|x)|x \bigr) \sqrt{nh^p
\alpha_n^{-1}} \bigl(\hat q_n(
\tau_j \alpha_n|x)-q(\tau_j
\alpha_n|x) \bigr) \Bigr\}_{j=1,\ldots,J}
\]
is asymptotically Gaussian, centered, with covariance matrix
$
\|K\|_2^2/g(x) \Sigma(x)
$
where $\Sigma_{j,j'}(x)= (\tau_j\tau_{j'})^{-\gamma(x)}\tau^{-1}_{j\wedge j'}$ for $(j,j')\in\{1,\ldots, J\}^2$.
\end{Theo}
Let us remark that, in the particular case where $J=1$, $\tau_1=1$ and
$\alpha_n=\alpha$ is
fixed in $(0,1)$, we find back the result of~\cite{BER}, Theorem 6.4.
Theorem~\ref{thquant} can be equivalently rewritten as
%
%co1 #&#
\begin{Coro}
\label{coroq}
Under the assumptions of Theorem~\ref{thquant}, the random vector
\[
\biggl\{ \sqrt{nh^p\alpha_n} \frac{q(\alpha_n|x)}{a(q(\alpha_n|x)|x)} \biggl(
\frac{\hat q_n(\tau_j
\alpha_n|x)}{q(\tau_j \alpha_n|x)}-1 \biggr) \biggr\}_{j=1,\ldots,J}
\]
is asymptotically Gaussian, centered, with covariance matrix
$
\|K\|_2^2/g(x) \tilde\Sigma(x)
$
where $\tilde\Sigma_{j,j'}(x) = (\tau_j\tau_{j'})^{-(\gamma
(x)\wedge
0)} \tau^{-1}_{j\wedge j'}$ for $(j,j')\in\{1,\ldots, J\}^2$.
\end{Coro}
Moreover,~\cite{deHaanFer}, Theorem 1.2.5 and~\cite{deHaanFer}, page~33, show that
%
%e5 #&#
\begin{equation}
\label{clef} \lim_{y\uparrow y_F(x)} \frac{a(y|x)}{y}= \gamma (x)\vee0.
\end{equation}
Under the assumptions of Theorem~\ref{thquant}, and from~(\ref{clef}),
it follows that
${\hat q_n(\tau_j \alpha_n|x)}/\break {q(\tau_j\alpha_n|x)} \toP1$ when
$n\to
\infty$ which can be
read as a weak consistency result for the considered estimator.
Besides, if $\gamma(x)>0$,
then collecting~(\ref{clef}) and Corollary~\ref{coroq} shows that the
random vector
\[
\biggl\{ \sqrt{nh^p\alpha_n} \biggl(\frac{\hat q_n(\tau_j \alpha_n|x)}{q(\tau_j \alpha_n|x)}-1
\biggr) \biggr\}_{j=1,\ldots,J}
\]
is asymptotically Gaussian, centered, with covariance matrix
$
\|K\|_2^2\gamma^2(x)/g(x) \tilde\Sigma(x)
$
where the coefficients of the covariance matrix can be simplified\vspace*{1.5pt}
$\tilde\Sigma_{j,j'}(x) = \tau^{-1}_{j\wedge j'}$ for $(j,j')\in\{
1,\ldots, J\}^2$.
Our results thus
build on and complement the analysis given by \cite{Test}, Theorem~2,
in the case
$\gamma(x)>0$.

As pointed out in~\cite{Test}, the condition $nh^p\alpha_n\to
\infty$ implies $\alpha_n> \log^p(n)/n$ eventually.
This condition provides a lower bound on the order of the extreme
conditional quantiles for the asymptotic normality of kernel estimators
to hold.
We now propose a scheme to estimate extreme
conditional quantiles without this restriction.
Let $\alpha_n\to0$ and $\beta_n/\alpha_n\to0$ as $n\to\infty$.
Suppose one has $\hat\gamma_n(x)$ and $\hat a_n(x)$ two estimators of
$\gamma(x)$ and $a(q(\alpha_n|x)|x)$, respectively.
Then, starting from the estimator $\hat q_n(\alpha_n|x)$ of $q(\alpha_n|x)$
defined in~(\ref{defestquant})
and making use of \eqref{F1quant},
it is possible to build
an estimator $\tilde q_n(\beta_n|x)$ of $q(\beta_n|x)$ which is an
extreme conditional quantile of higher order than $q(\alpha_n|x)$:
%
%e6 #&#
\begin{equation}
\label{qtilde} \tilde q_n(\beta_n|x)= \hat
q_n(\alpha_n|x) + K_{\hat\gamma
_n(x)}(\alpha_n/
\beta_n) \hat a_n(x).
\end{equation}
Let us consider, for all $z\in\R$, the function defined for all $u>1$ by
\[
K'_z(u)=\frac{\partial K_z(u)}{\partial z} = \int
_1^u v^{z-1} \log(v) \,\mathrm{d}v.
\]
The following result provides a quantile regression analog of \cite{deHaanFer}, Theorem~4.3.1.
%
%th2 #&#
\begin{Theo}
\label{thquant2}
Suppose \textup{(A.1)} holds and let $\alpha_n\to0$, $\beta_n/\alpha_n\to
0$. Let
$\hat q_n( \alpha_n|x)$
be the kernel estimator of $q(\alpha_n|x)$ defined in~(\ref{defestquant}).
Let $\hat\gamma_n(x)$ and $\hat a_n(x)$ be two estimators of $\gamma
(x)$ and $a(q(\alpha_n|x)|x)$, respectively,
such that
%
%e7 #&#
\begin{equation}
\label{condtriplet} \Lambda_n^{-1} \biggl( \hat
\gamma_n(x)-\gamma(x), \frac{\hat a_n(x)}{a(q(\alpha_n|x)|x)}-1, \frac{\hat q_n( \alpha_n|x)-q( \alpha_n|x)}{a(q(\alpha_n|x)|x)}
\biggr)^t \tod\zeta(x),
\end{equation}
where $\zeta(x)$ is a nondegenerate
$\R^3$ random vector,
\[
\Lambda_n \log(\alpha_n/\beta_n)\to0 \quad\mbox{and}\quad \Lambda_n^{-1} \frac
{b(\beta_n/\alpha_n,\alpha_n|x)}{K'_{\gamma(x)}(\alpha_n/\beta_n)}\to0
\]
as $n\to\infty$. Then,
\[
\Lambda_n^{-1} \biggl(\frac{\tilde q_n(\beta_n|x)-q(\beta_n|x)}{a(q(\alpha_n|x)|x) K{}_{\gamma(x)}'(\alpha_n/\beta_n)} \biggr) \tod
c(x)^t \zeta(x),
\]
where $c(x)^t = (1, - (\gamma(x)\wedge0), (\gamma(x)\wedge0)^2
)$.
\end{Theo}
As an illustration, for all $r\in(0,1)$, let us consider
$\tau_j=r^{j-1}$, $j=1,\ldots,J$.
The following estimators of $\gamma(x)$ and $a(q(\alpha_n|x)|x)$ are
introduced
\begin{eqnarray*}
\hat{\gamma}_n^{\mathrm{RP}}(x)&=& \frac{1}{\log r}\sum
_{j=1}^{J-2} \pi_j \log \biggl(
\frac{\hat q_n(\tau_j \alpha_n|x) - \hat
q_n(\tau_{j+1}\alpha_n|x)} {
\hat q_n(\tau_{j+1} \alpha_n|x) - \hat q_n(\tau_{j+2}\alpha_n|x)} \biggr),
\\
\hat{a}_n^{\mathrm{RP}}(x)&=& \frac{1}{K_{\hat{\gamma}_n^{\mathrm
{RP}}(x)}(r)} \sum
_{j=1}^{J-2} \pi_j r^{
\hat{\gamma}_n^{\mathrm{RP}}(x)j} \bigl(
\hat q_n( \tau_j \alpha_n|x) - \hat
q_n( \tau_{j+1}\alpha_n|x) \bigr),
\end{eqnarray*}
where $(\pi_j)$ is a sequence of weights summing to one. Let us
highlight that
$\hat{\gamma}_n^{\mathrm{RP}}(x)$ is an adaptation to the
conditional case of the Refined
Pickands estimator introduced in~\cite{Drees}.
The joint asymptotic normality of $(\hat{\gamma}_n^{\mathrm
{RP}}(x),\hat{a}_n^{\mathrm{RP}}(x),\hat q_n(
\alpha_n|x))$
is established in the next theorem.

%th3 #&#
\begin{Theo}
\label{triplet}
Suppose \textup{(A.1)}, \textup{(A.2)} and \textup{(A.3)} hold. Let $x\in\R^p$ such that $g(x)>0$.
If $\alpha_n \to0$ and there exists $\kappa\in(0,\tau_J)$ such that
\[
nh^p\alpha_n \to\infty,  nh^{p}
\alpha_n \Biggl(h \vee\Delta_{\kappa
}(x,\alpha_n)
\vee\bigvee_{j=1}^J b(\tau_j,
\alpha_n|x) \Biggr)^2 \to0
\]
as $n\to\infty$, then the random vector
\[
\sqrt{nh^{p}\alpha_n } \biggl( \hat{\gamma}_n^{\mathrm
{RP}}(x)-
\gamma(x), \frac{\hat{a}_n^{\mathrm{RP}}(x)}{a(q(\alpha_n|x)|x)}-1, \frac
{\hat q_n( \alpha_n|x)-q( \alpha_n|x)}{a(q(\alpha_n|x)|x)} \biggr)^t
\]
is asymptotically centered and Gaussian.
\end{Theo}
The asymptotic covariance matrix is denoted by $S(x)$. It can be explicitly
calculated from~(\ref{matS}) in the proof of Theorem~\ref{triplet}, but
the result would be too
complicated to be reported here.
As a consequence of the two above theorems, one obtains the asymptotic
normality of
the extreme conditional quantile estimator built on $\hat{\gamma
}_n^{\mathrm{RP}}(x)$ and
$\hat{a}_n^{\mathrm{RP}}(x)$:
\[
\tilde{q}_n^{\mathrm{RP}}(\beta_n|x):= \hat
q_n(\alpha_n|x) + K_{\hat{\gamma}_n^{\mathrm{RP}}(x)}( \alpha_n/
\beta_n) \hat{a}_n^{\mathrm{RP}}(x).
\]
%
%co2 #&#
\begin{Coro}
Suppose \textup{(A.1)}, \textup{(A.2)} and \textup{(A.3)} hold. Let $x\in\R^p$ such that $g(x)>0$.
If $\alpha_n \to0$, $\beta_n/\alpha_n\to0$ and there exists
$\kappa\in
(0,\tau_J)$ such that
\[
\frac{nh^p\alpha_n}{(\log(\alpha_n/\beta_n))^2} \to\infty, \qquad nh^{p}\alpha_n \Biggl(h
\vee\Delta_{\kappa}(x,\alpha_n) \vee\bigvee
_{j=1}^J b(\tau_j,\alpha_n|x)
\vee\frac{b(\beta_n/\alpha_n,\alpha_n|x)}{K'_{\gamma(x)}(\alpha_n/\beta_n)} \Biggr)^2 \to0
\]
as $n\to\infty$, then
\[
\sqrt{nh^p\alpha_n} \biggl(\frac{\tilde{q}_n^{\mathrm{RP}}(\beta_n|x)-q(\beta_n|x)}{a(q(\alpha_n|x)|x)
K'_{\gamma(x)}(\alpha_n/\beta_n)} \biggr)
\]
is asymptotically Gaussian, centered with variance $c(x)^t S(x) c(x)$.
\end{Coro}

Finally, two particular cases of $\hat{\gamma}_n^{\mathrm{RP}}(x)$
may be
considered. First, constant weights $\pi_1=\cdots=\pi_{J-2}=1/(J-2)$
yield
\[
\hat{\gamma}_n^{\mathrm{RP},1}(x)= \frac{1}{(J-2)\log r} \log \biggl(
\frac{\hat q_n(\tau_1 \alpha_n|x) - \hat q_n(\tau_{2}\alpha_n|x)} {
\hat q_n(\tau_{J-1} \alpha_n|x) - \hat q_n(\tau_{J}\alpha_n|x)} \biggr).
\]
Clearly, when $J=3$, this estimator reduces the kernel Pickands
estimator introduced and studied
in~\cite{Test} in the situation where $\gamma(x)>0$. Second,
linear\vadjust{\goodbreak}
weights $\pi_j=2j/((J-1)(J-2))$
for $j=1,\ldots,J-2$ give rise to a new estimator
\[
\hat{\gamma}_n^{\mathrm{RP},2}(x)= \frac{2}{(J-1)(J-2)\log r}\sum
_{j=1}^{J-2} \log \biggl( \frac{\hat q_n(\tau_j \alpha_n|x) - \hat q_n(\tau_{j+1}\alpha_n|x)} {
\hat q_n(\tau_{J-1} \alpha_n|x) - \hat q_n(\tau_{J}\alpha_n|x)} \biggr),
\]
which can be read as the average of $J-1$ estimators $\hat{\gamma
}_n^{\mathrm{RP},1}(x)$.
These estimators are now compared on finite sample situations.

%--------------------------------------------------------------------

%s4 #&#
\section{Some simulation evidence}\label{mc}
Section~\ref{MCE} provides Monte Carlo evidence that the extreme
quantile function
estimator $\tilde{q}_n^{\mathrm{RP},1}(\beta_n|x)$ is efficient
relative to the
version $\tilde{q}_n^{\mathrm{RP},2}(\beta_n|x)$, whether $\gamma
(x)$ is
positive, negative or zero, and outperforms the estimator $\hat
{q}_n(\beta_n|x)$ for heavy-tailed conditional distributions.
Section~\ref{BG} provides a comparison with the promising local
smoothing
%polynomial maximum likelihood
approach introduced in~\cite{BEG} and \cite{BES}, Section~7.5.2.
Practical guidelines for selecting the bandwidth $h$ and the order
$\alpha_n$ are suggested in Section~\ref{PG}.

%s4.1 #&#
\subsection{Monte Carlo experiments}
\label{MCE}
To evaluate finite-sample performance of the conditional extreme-value
index and
extreme quantile estimators described above, we have undertaken some simulation
experiments following the model
\[
Y_i = \mathcal{G}(X_i) + \sigma(X_i)
U_i, \qquad i=1,\ldots,n.
\]
The local scale factor, $\sigma(x)=(1+x)/10$, is linearly increasing
in $x$,
while the local location parameter
\[
\mathcal{G}(x)=\sqrt{x(1-x)} \sin \biggl(\frac{2\uppi
(1+2^{-7/5})}{x+2^{-7/5}} \biggr)
\]
has been introduced in~\cite{RUP}, Section~17.5.1.
The design points $X_i$ are generated following a standard uniform distribution.
The $U_i$'s are independent and their conditional distribution given
$X_i=x$ is chosen to be standard Gaussian, Student $t_{k(x)}$, or
Beta$(\nu(x),\nu(x))$, with
\[
k(x)= \bigl[\nu(x) \bigr]+1, \qquad \nu(x)= \bigl\{ \bigl(\tfrac
{1}{10}+\sin(\uppi
x) \bigr) \bigl(\tfrac
{11}{10}-\tfrac{1}{2}\exp \bigl
\{-64(x-1/2)^2 \bigr\} \bigr) \bigr\}^{-1},
\]
and $[\nu(x)]$ being the integer part of $\nu(x)$. Let us recall
that the Gaussian distribution belongs to the Gumbel MDA, {that
is}, $\gamma(x)=0$, the Student distribution $t_{k(x)}$
belongs to the Fr\'echet MDA with $\gamma(x)=1/k(x)>0$ and the Beta
distribution belongs to the Weibull MDA with
$\gamma(x)=-1/\nu(x)<0$.

In all cases, we have
$q(\beta|x)=\mathcal{G}(x)+\sigma(x)\bar{F}^{\leftarrow
}_{U|X}(\beta
|x)$, for $\beta\in(0,1)$.
All the experiments were performed over $400$ simulations for\vadjust{\goodbreak} $n=200$,
and the kernel function $K$ was chosen to be the Triweight kernel
\[
K(t)=\tfrac{35}{32} \bigl(1-t^2 \bigr)^{3}\indic\{-1
\leq t\leq1\}.
\]
Monte Carlo experiments were first devoted to accuracy of the two conditional
extreme-value index estimators $\hat{\gamma}_n^{\mathrm{RP},1}(x)$
and $\hat{\gamma}_n^{\mathrm{RP},2}(x)$.
The measures of efficiency for each simulation used were the mean
squared error and the bias
\[
\operatorname{MSE} \bigl\{\hat\gamma_n(\cdot) \bigr\}=
\frac{1}{L}\sum_{\ell=1}^L \bigl\{\hat
\gamma_n(x_\ell)-\gamma(x_\ell) \bigr
\}^2,\qquad  \operatorname{Bias} \bigl\{\hat\gamma_n(\cdot) \bigr
\}= \frac{1}{L}\sum_{\ell=1}^L \bigl\{
\hat\gamma_n(x_\ell)-\gamma(x_\ell) \bigr\}
\]
for $\hat\gamma_n(x)=\hat{\gamma}_n^{\mathrm{RP},1}(x), \hat
{\gamma}_n^{\mathrm{RP},2}(x)$, with the
$x_\ell$'s being $L=100$ points regularly distributed in $[0,1]$.
To guarantee a fair comparison among the two estimation methods, we
used for
each estimator the parameters $(\alpha_n,h)$ minimizing its mean
squared error,
with $\alpha_n$ ranging over $\mathcal{A}=\{0.1,0.15,0.2,\ldots
,0.95\}$
and the bandwidth $h$ ranging over a grid $\mathcal{H}$ of $50$ points
regularly
distributed between $h_{\mathrm{min}}=\max_{1\leq i<n}|X_{(i+1)}-X_{(i)}|$ and
$h_{\mathrm{max}}=|X_{(n)}-X_{(1)}|/2$, where $X_{(1)}\leq\cdots\leq X_{(n)}$
are the
ordered observations.
The resulting values of MSE and bias are averaged on the $400$ Monte
Carlo replications
and reported in Table~\ref{tabb1} for
$J\in\{3,4,5\}$ and $r\in\{1/J,(J-1)/J\}$.

%t1 #&#
\begin{table}
\def\arraystretch{0.9}
\caption{Performance of $\hat{\gamma}_n^{\mathrm{RP},1}(x)$ and
$\hat{\gamma}_n^{\mathrm{RP},2}(x)$ -- Results
averaged on $400$ simulations with $n= 200$.
The results may not be available for $r=1/J$
and $J=5$ since the numerator $\{\hat{q}_n(\tau_j\alpha_n|x)-\hat
{q}_n(\tau_{j+1}\alpha_n|x)\}$ and the denominator $\{\hat{q}_n(\tau_{J-1}\alpha_n|x)-\hat{q}_n(\tau_{J}\alpha_n|x)\}$ in the
definitions of both estimators might be null when $n$ is not large enough}
\label{tabb1}
\begin{tabular*}{\textwidth}{@{\extracolsep{\fill
}}ld{1.4}d{1.4}d{2.4}d{2.4}@{}}
\hline
& \multicolumn{2}{l}{MSE}& \multicolumn{2}{l@{}}{Bias}\\[-4pt]
& \multicolumn{2}{c}{\hrulefill}& \multicolumn{2}{c@{}}{\hrulefill
}\\
& \multicolumn{1}{c}{$\hat{\gamma}_n^{\mathrm{RP},1}(x)$} &
\multicolumn{1}{c}{$\hat{\gamma}_n^{\mathrm{RP},2}(x)$}
& \multicolumn{1}{c}{$\hat{\gamma}_n^{\mathrm{RP},1}(x)$} &
\multicolumn{1}{c@{}}{$\hat{\gamma}_n^{\mathrm{RP},2}(x)$} \\
\hline
& \multicolumn{4}{c}{$ r=1/J$}\\
{Gaussian} & & & &\\
$J=3$ & 0.2026 & 0.2026 & -0.2415 & -0.2415 \\
$J=4$ & 0.1915 & 0.2018 & -0.3270 & -0.3501 \\
$J=5$ & \multicolumn{1}{c}{NaN} & \multicolumn{1}{c}{NaN} &
\multicolumn{1}{c}{NaN} & \multicolumn{1}{c}{NaN} \\[3pt]
{Student} & & & &\\
$J=3$ & 0.2882 & 0.2882 & -0.2964 & -0.2964 \\
$J=4$ & 0.3350 & 0.2837 & -0.4167 & -0.3480 \\
$J=5$ & \multicolumn{1}{c}{NaN} & \multicolumn{1}{c}{NaN} &
\multicolumn{1}{c}{NaN} & \multicolumn{1}{c}{NaN} \\[3pt]
{Beta} & & & &\\
$J=3$ & 0.1157 & 0.1157 & -0.0730 & -0.0730 \\
$J=4$ & 0.0510 & 0.0597 & -0.0811 & -0.0750 \\
$J=5$ & \multicolumn{1}{c}{NaN} & \multicolumn{1}{c}{NaN} &
\multicolumn{1}{c}{NaN} & \multicolumn{1}{c}{NaN}
\\[6pt]
& \multicolumn{4}{c}{$r=(J-1)/J$}\\
{Gaussian} & & & &\\
$J=3$ & 0.7656 & 0.7656 & -0.3213 & -0.3213 \\
$J=4$ &0.6730 & 0.7960 & -0.3455 & -0.3747 \\
$J=5$ &0.7305 & 0.9128 & -0.4104 & -0.4107 \\[3pt]
{Student} & & & &\\
$J=3$ &1.1109 & 1.1109 & -0.4497 & -0.4497 \\
$J=4$ &0.9991 & 1.1997 & -0.4384 & -0.4597 \\
$J=5$ &1.1245 & 1.3331 & -0.5715 & -0.5872 \\[3pt]
{Beta} & & & &\\
$J=3$ &0.6737 & 0.6737 & -0.2591 & -0.2591 \\
$J=4$ &0.5861 & 0.6891 & -0.2338 & -0.2432 \\
$J=5$ & 0.6431 & 0.8167 & -0.2185 & -0.2757 \\
\hline
\end{tabular*}
%
%%%%%%%%%%%%%%%%%%%%%
\end{table}

It does appear that the results for $r=1/J$ are superior to those
for $r=(J-1)/J$, uniformly in $J$. For these desirable results,
it may be seen that the estimator $\hat{\gamma}_n^{\mathrm
{RP},1}(x)$ performs better than
$\hat{\gamma}_n^{\mathrm{RP},2}(x)$ in the Gaussian error model,
whereas the latter is
superior to the former in the Student error model. It may be also seen that
there is no winner in the Beta error model in terms of both MSE and bias.

Turning to the performance of the extreme conditional quantile
estimators, we
consider as above the two measures of performance
\begin{eqnarray*}
\operatorname{MSE} \bigl\{q_n(\beta_n|\cdot) \bigr\}&=&
\frac{1}{L}\sum_{\ell=1}^L \bigl
\{q_n(\beta_n|x_\ell)-q(\beta_n|x_\ell
) \bigr\}^2 ,\\
 \operatorname{Bias} \bigl\{q_n(
\beta_n|\cdot) \bigr\}&=& \frac{1}{L}\sum
_{\ell
=1}^L \bigl\{ q_n(
\beta_n|x_\ell)-q(\beta_n|x_\ell)
\bigr\}
\end{eqnarray*}
for $q_n(\beta_n|x)=\hat{q}_n(\beta_n|x)$, $\tilde{q}_n^{\mathrm
{RP},1}(\beta_n|x)$, $\tilde{q}_n^{\mathrm{RP},2}(\beta_n|x)$.
The averaged MSE and bias of these three estimators of $q(\beta_n|x)$,
computed for $\beta_n\in\{0.05,0.01,0.005\}$, $J\in\{3,4\}$ and
$r = 1/J$, over $400$ Monte Carlo simulations are displayed in
Table~\ref{tabb2}. Here also, we used for each estimator the smoothing
parameters $(\alpha_n,h)$ minimizing its MSE over the grid of values
$\mathcal{A}\times\mathcal{H}$ described above.

When comparing the estimators $\tilde{q}_n^{\mathrm{RP},1}(\beta_n|x)$ and
$\tilde{q}_n^{\mathrm{RP},2}(\beta_n|x)$ themselves with
$\hat{q}_n(\beta_n|x)$, the results (both in terms of MSE and bias)
indicate that
$\tilde{q}_n^{\mathrm{RP},2}(\beta_n|x)$ is slightly less efficient than
$\tilde{q}_n^{\mathrm{RP},1}(\beta_n|x)$ in all cases, and that the
latter is appreciably more efficient
than $\hat{q}_n(\beta_n|x)$ only in the Student error model.
It may be also noticed that $\hat{q}_n(\beta_n|x)$ is more
efficient but not by much (especially when $J=3$) in the Gaussian and
Beta error models.

%However, both $\qRPu$ and
%$\qRPd$ appear to be appreciably more efficient
%than $\hat{q}_n(\beta_n|x)$ (especially for $r=1/J$ and $J=3$) in the
%Student error model.

%
% In terms of bias, it may be seen in the Gaussian error model that
% $\tilde{q}^{RP,1}_n(\beta_n|\cdot)$ outperforms both
% $\tilde{q}^{RP,2}_n(\beta_n|\cdot)$ and $\hat{q}_n(\beta_n|\cdot)$
%for some
% values of $r$ and $J$, especially $r=(J-1)/J$ and $J=4$.
% In contrast, in the Student error model, there is no winner between
% $\tilde{q}^{RP,1}_n(\beta_n|\cdot)$ and $\tilde{q}^{RP,2}_n(\beta_n|
% but both are superior to $\hat{q}_n(\beta_n|\cdot)$ almost overall.
% Finally, in the Beta error model, $\tilde{q}^{RP,1}_n(\beta_n|\cdot)$
%and
% $\tilde{q}^{RP,2}_n(\beta_n|\cdot)$ seem to achieve nice bias
%estimates at the
% specific values $r=(J-1)/J$ and $J=3$, although $\hat{q}_n(\beta_n|
% affords slightly better estimates.
%We obtain the same conclusions for the higher order $
%Table \ref{tabzz}, for both bias and MSE.

%t2 #&#
\begin{table}
\def\arraystretch{0.9}
\caption{Performance of
$\tilde{q}_n^{\mathrm{RP},1}(\beta_n|x), \tilde{q}_n^{\mathrm
{RP},2}(\beta_n|x)$
and $\hat{q}_n(\beta_n|x)$ with $\beta_n=0.05$ (top), $\beta_n=0.01$
(middle) and $\beta_n=0.005$ (bottom) -- Results averaged on $400$
simulations with $n= 200$}
\label{tabb2}
\begin{tabular*}{\textwidth}{@{\extracolsep{\fill}}lllld{2.4}d{2.4}d{2.4}@{}}
\hline
& \multicolumn{3}{l}{MSE}
& \multicolumn{3}{l@{}}{Bias}\\[-4pt]
& \multicolumn{3}{c}{\hrulefill}
& \multicolumn{3}{c@{}}{\hrulefill}\\
& \multicolumn{1}{c}{$\tilde{q}_n^{\mathrm{RP},1}(\beta_n|x)$} &
\multicolumn{1}{c}{$\tilde{q}_n^{\mathrm{RP},2}(\beta_n|x)$} &
\multicolumn{1}{c}{$\hat{q}_n(\beta_n|x)$}
& \multicolumn{1}{c}{$\tilde{q}_n^{\mathrm{RP},1}(\beta_n|x)$} &
\multicolumn{1}{c}{$\tilde{q}_n^{\mathrm{RP},2}(\beta_n|x)$} &
\multicolumn{1}{c@{}}{$\hat{q}_n(\beta_n|x)$} \\
\hline
& \multicolumn{6}{c}{$r=1/J$, $\beta_n=0.05$}\\
{Gaussian} & & & & & & \\
$J=3$ & 0.0110 & 0.0110 & 0.0108 & 0.0001 & 0.0001 & 0.0063 \\
$J=4$ & 0.0591 & 0.0796 & 0.0108 & 0.1136 & 0.1131 & 0.0063 \\[3pt]
{Student} & & & & & &\\
$J=3$ & 0.0307 & 0.0307 & 0.0771 & -0.0134 & -0.0134 & 0.0871 \\
$J=4$ & 0.0532 & 0.0743 & 0.0771 & 0.0792 & 0.0792 & 0.0871 \\[3pt]
{Beta} & & & & & &\\
$J=3$ & 0.0091 & 0.0091 & 0.0022 & 0.0505 & 0.0505 & 0.0135 \\
$J=4$ & 0.0745 & 0.1002 & 0.0022 & 0.1746 & 0.1752 & 0.0135 \\[6pt]
& \multicolumn{6}{c}{$r=1/J$, $\beta_n=0.01$} \\
{Gaussian} & & & & & & \\
$J=3$ & 0.0265 & 0.0265 & 0.0161 & -0.0776 & -0.0776 & -0.0360 \\
$J=4$ & 0.0693 & 0.0926 & 0.0161 & 0.1092 & 0.1225 & -0.0360 \\[3pt]
{Student} & & & & & &\\
$J=3$ & 0.1115 & 0.1115 & 0.6825 & -0.0895 & -0.0895 & -0.0959 \\
$J=4$ & 0.1304 & 0.3992 & 0.6825 & 0.0018 & 0.1089 & -0.0959 \\[3pt]
{Beta} & & & & & &\\
$J=3$ & 0.0143 & 0.0143 & 0.0034 & 0.0523 & 0.0523 & 0.0212 \\
$J=4$ & 0.1038 & 0.1265 & 0.0034 & 0.1964 & 0.2064 & 0.0212 \\[6pt]
%%%%%%%%%%%%%%%%%%%%%%%%%%%%%%%%%%%%%%%%%%%%%%
& \multicolumn{6}{c}{$r=1/J$, $\beta_n=0.005$}\\
{Gaussian} & & & & & & \\
$J=3$ & 0.0354 & 0.0354 & 0.0203 & -0.0981 & -0.0981 & -0.0524 \\
$J=4$ & 0.0719 & 0.0932 & 0.0203 & 0.0982 & 0.1073 & -0.0524 \\[3pt]
{Student} & & & & & &\\
$J=3$ & 0.2919 & 0.2919 & 0.9782 & -0.1623 & -0.1623 & -0.2605\\
$J=4$ & 0.4569 & 0.9748 & 0.9782 & -0.1920 & 0.0280 & -0.2605\\[3pt]
{Beta} & & & & & &\\
$J=3$ & 0.0155 & 0.0155 & 0.0038 & 0.0536 & 0.0536 &0.0239\\
$J=4$ & 0.1130 & 0.1337 & 0.0038 & 0.1871 & 0.2111 &0.0239\\
\hline
%%%%%%%%%%%%%%%%%%%%
\end{tabular*}
\end{table}

%----------------------------------------------------

%s4.2 #&#
\subsection{\texorpdfstring{Benchmark nonparametric estimators of $\gamma(x)$ and $q(\beta_n|x)$}
{Benchmark nonparametric estimators of gamma(x) and q(beta n|x)}}\label{BG}

Alternative modern smoothing techniques were discussed in,
for example, \cite{BES}, Section~7.5. For comparison, we focus on the
prominent
local polynomial maximum likelihood estimation. This contribution fits a
generalized Pareto (GP) model to the exceedances $Z^x_i=Y_j-u_x$ given
$Y_j>u_x$, for a high threshold $u_x$, where $j$ denotes the original
index of
the $i$th exceedance. Let $N_x$ be the number of all exceedances over $u_x$
and rearrange the indices of the explanatory variable such that $X_i$ denotes
the covariate observation associated with exceedance $Z^x_i$. If
$\tilde{g}(z;\sigma,\gamma)$ stands for the GP density, then the local
polynomial
maximum likelihood approach maximizes the kernel weighted log-likelihood
function
\[
L_{N_x}(\beta_1,\beta_2)= \frac{1}{N_x}
\sum_{i=1}^{N_x}\log\tilde{g}
\Biggl(Z^x_i;\sum_{j=0}^{p_1}
\beta_{1j}(X_i-x)^j,\sum
_{j=0}^{p_2}\beta_{2j}(X_i-x)^j
\Biggr)K_h(X_i-x)
\]
with respect to $(\beta_1',\beta_2')=(\beta_{10},\ldots,\beta_{1p_1},\beta_{20},\ldots,\beta_{2p_2})$ to get the estimates
$\hat{\sigma}^{\mathrm{GP}}_n(x)=\hat{\beta}_{10}$ and
$\hat{\gamma}^{\mathrm{GP}}_n(x)=\hat{\beta}_{20}$ of the
parameter functions
$\sigma(x)$
and $\gamma(x)$ of the GP distribution fitted to the
exceedances over $u_x$. Note that local polynomial fitting also provides
estimates of the derivatives of $\sigma(x)$ and $\gamma(x)$ up to
order $p_1$
and $p_2$, respectively.
%Given that our interest is only in estimation of the
%parameter functions, we confine ourselves to $p_1=p_2=0$.
%%%%%%%%%%
In order to not overload the estimation
procedure, we confine ourselves to $p_1=p_2=0$.
%%%%%%%%%%
The Monte Carlo results for
$\hat{\gamma}^{\mathrm{GP}}_n(\cdot)$ are reported in Table~\ref
{Tablerev} (top).
For each simulation, we used the parameters $(h,u)$ that minimize the
$\operatorname{MSE}\{\hat{\gamma}^{\mathrm{GP}}_n(\cdot)\}$, with
the bandwidth ranging
over the
grid $\mathcal{H}$ described above and the threshold ranging over the
$\alpha$th sample quantiles of $Y$, where $\alpha\in\mathcal{A}$.
%%%%%%%%
%%%%%%%%
The estimator $\hat{\gamma}^{\mathrm{GP}}_n$ has clearly smaller
MSEs than the
$\hat{\gamma}^{\mathrm{RP}}_n$ estimators in the Gaussian and
Student error
models, but
it seems to be less efficient in the Beta error model than both
$\hat{\gamma}^{\mathrm{RP},1}_n$ and $\hat{\gamma}^{\mathrm
{RP},2}_n$ for $J=4$ and $r=1/J$.
From a theoretical point of view, it should be clear that the pointwise
asymptotic normality of $\hat{\gamma}^{\mathrm{GP}}_n(x)$ is proved
in \cite{BEG} only
in case $\gamma(x)>0$. Moreover, the proof is restricted to the
setting where
the design points $X_i$ are deterministic.

%t3 #&#
\begin{table}
\caption{Performance of $\hat{\gamma}^{\mathrm{GP}}_n$ and
$\hat
{q}^{\mathrm{GP}}_n(\beta_n|\cdot)$ -- Results averaged on $400$ simulations
with $n= 200$}\label{Tablerev}
\begin{tabular*}{200pt}{@{\extracolsep{\fill}}lll@{}}
\hline
&\multicolumn{1}{l}{$\operatorname{MSE}\{\hat{\gamma}^{\mathrm
{\mathrm{GP}}}_n\}$}
&
\multicolumn{1}{l@{}}{$\operatorname{Bias}\{\hat{\gamma}^{\mathrm
{\mathrm{GP}}}_n\}$}\\
\hline
{Gaussian} & 0.1324 & $-0.2671$ \\
{Student} & 0.1310 & $-0.2238$ \\
{Beta} & 0.0675 & $-0.0221$ \\
\hline
\end{tabular*}
\vspace*{6pt}
\begin{tabular*}{\textwidth}{@{\extracolsep{\fill}}lllllll@{}}
\hline
&\multicolumn{2}{l}{$\beta_n=0.05$}
&
\multicolumn{2}{l}{$\beta_n=0.01$}
&
\multicolumn{2}{l@{}}{$\beta_n=0.005$}\\[-4pt]
&\multicolumn{2}{c}{\hrulefill}
&
\multicolumn{2}{c}{\hrulefill}
&
\multicolumn{2}{c@{}}{\hrulefill}\\
&
\multicolumn{1}{l}{$\operatorname{MSE}$}
&
\multicolumn{1}{l}{$\operatorname{Bias}$}
&
\multicolumn{1}{l}{$\operatorname{MSE}$}
&
\multicolumn{1}{l}{$\operatorname{Bias}$}
&
\multicolumn{1}{l}{$\operatorname{MSE}$}
&
\multicolumn{1}{l@{}}{$\operatorname{Bias}$}
\\
\hline
{Gaussian} & 0.0184 & 0.0974 & 0.0278 & 0.0952 & 0.0315 & \phantom
{$-$}0.0861 \\
{Student} & 0.1346 & 0.1526 & 0.6924 & 0.0895 & 1.0232 & $-0.0452$ \\
{Beta} & 0.0364 & 0.1578 & 0.0659 & 0.2067 & 0.0786 & \phantom
{$-$}0.2242 \\
\hline
\end{tabular*}
\end{table}

%%%%%%%%%%%%%%%%%%%%%%%%%%

\begin{figure}

\includegraphics{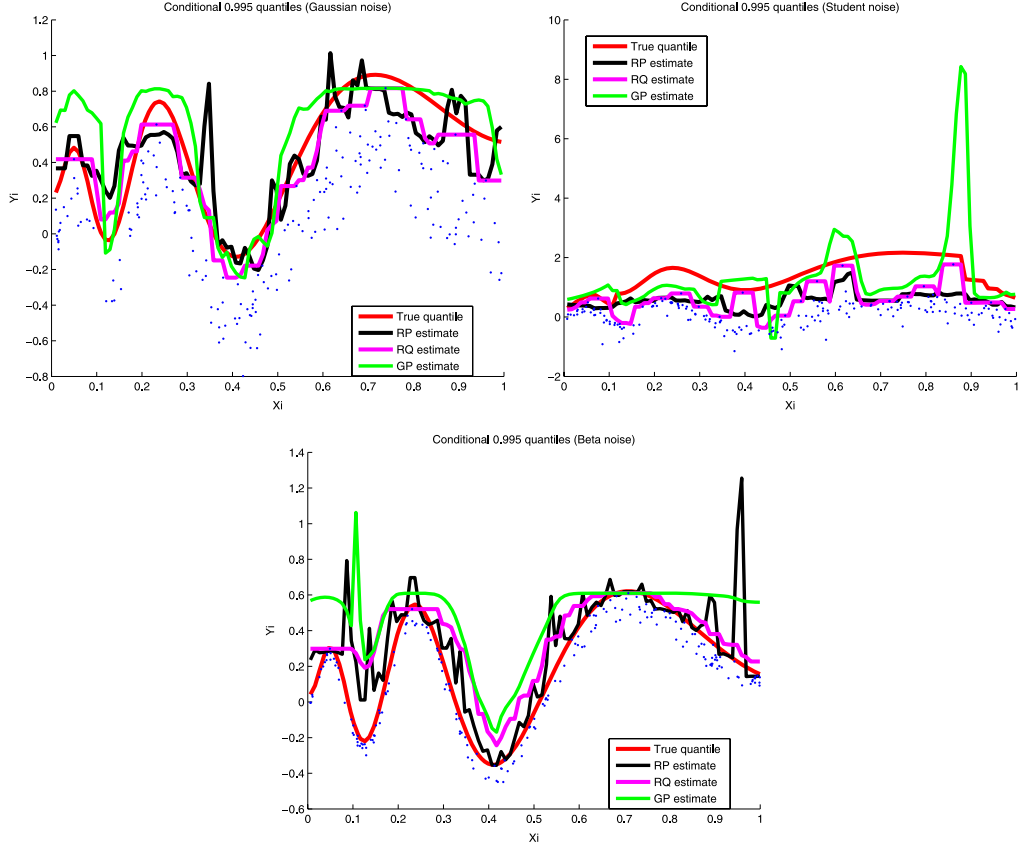}

\caption{Typical realizations for simulated samples of size $n=200$. From
left to right and from top to bottom, $Y|X$ is Gaussian, Student,
Beta. The true quantile function $q(\beta_n|\cdot)$ in red with
$\beta_n=1/n$. Its estimators $\hat{q}_n(\beta_n|\cdot)$ in magenta,
$\tilde{q}^{\mathrm{RP},1}_n(\beta_n|\cdot)\equiv\tilde
{q}^{\mathrm{RP},2}_n(\beta_n|\cdot)$
in black with $r=1/J$ and $J=3$, and $\hat{q}^{\mathrm{GP}}_n(\beta_n|\cdot)$ in
green. The observations $(X_i,Y_i)$ are depicted as blue points.}
\label{figA}
\end{figure}

On the other hand, as suggested in~\cite{BES}, Section~7.5.2 and \cite{BEG},
the extreme conditional quantile $q(\beta_n|x)$ can be estimated by
\[
\hat{q}^{\mathrm{GP}}_n(\beta_n|x) := u_x
+ \frac{\hat{\sigma
}^{\mathrm{GP}}_n(x)}{\hat
{\gamma}^{\mathrm{GP}}_n(x)} \biggl[ \biggl(\frac{n^{\star
}_{xh}\beta_n}{k_x} \biggr)^{-\hat{\gamma}^{\mathrm{GP}}_n(x)}-1
\biggr],
\]
where $n^{\star}_{xh}$ is the number of observations in $[x-h,x+h]$
and $k_x$
is the number of exceedances receiving positive weight.
Table~\ref{Tablerev} (bottom) reports the Monte Carlo estimates
obtained by
using in each simulation the parameters $(h,u)$ that minimize the
$\operatorname{MSE}\{\hat{q}^{\mathrm{GP}}_n(\beta_n|\cdot)\}$,
where $h\in
\mathcal
{H}$ and
$u$ ranges over the $\alpha$th sample quantiles of $Y$ with $\alpha
\in
\mathcal{A}$.
%%%%%%%%%%
%%%%%%%%%%
In all cases, the regression quantile (RQ) estimator
$\hat{q}_n(\beta_n|\cdot)$ do appear to be more efficient than
$\hat{q}^{\mathrm{GP}}_n(\beta_n|\cdot)$. Compared with the
$\tilde{q}^{\mathrm{RP}}_n(\beta_n|\cdot)$ estimators (for $J=3$
and $r=1/J$),
$\hat{q}^{\mathrm{GP}}_n(\beta_n|\cdot)$ seems to be more efficient
only in the
Gaussian error model for $\beta_n=0.005=1/n$, but not by much.
A typical realization of the experiment in each simulated scenario is shown
in Figure~\ref{figA}, where the smoothing parameters of each
estimator were
chosen in such a way to minimize its MSE.
From a theoretical viewpoint, unlike our estimators $\hat{q}_n(\beta_n|x)$ and $\tilde{q}^{\mathrm{RP}}_n(\beta_n|x)$, the asymptotic
distribution of $\hat{q}^{\mathrm{GP}}_n(\beta_n|x)$ is not elucidated
yet.
%%%%%%%%%%%%%%%%%%%
%
%f1 #&#

%%%%%%%%%%%%%%%

%s4.3 #&#
\subsection{\texorpdfstring{Data-driven rules for selecting the parameters $h$ and $\alpha_n$}
{Data-driven rules for selecting the parameters h and alpha n}}\label{PG}
The use of the `RQ' estimator
$\hat{q}_n(\beta_n|x)\equiv\hat{q}_n(\beta_n|x;h)$, which relies
on the
inversion of $\hat{\hspace*{-2pt}\bar F}_n(\cdot|\cdot)$, requires only the choice
of the
bandwidth $h$ in an interval $\mathcal{H}$ of lower and upper bounds
given,
respectively, by, say,
$h_{\mathrm{min}}:=\max_{1\leq i<n} (X_{(i+1)}-X_{(i)} )$ and
$h_{\mathrm{max}}:= (X_{(n)}-X_{(1)} )/4$.
One way to select this parameter is by employing the cross-validation
criterion as in~\cite{Test} to obtain
\[
\label{cv} h_{\mathrm{cv}}=\mathop{\arg\min}_{h\in\mathcal{H}}\sum
_{i=1}^n\sum_{j=1}^n
\bigl\{ \indic(Y_i\geq Y_j)-\hat{\hspace*{-2pt}\bar
F}_{n,-i}(Y_j|X_i) \bigr\}^2,
\]
where $\hat{\hspace*{-2pt}\bar F}_{n,-i}(\cdot|\cdot)$ is the estimator $\hat{\hspace*{-2pt}\bar
F}_n(\cdot|\cdot)$ computed from the sample $\{(X_{j},Y_{j}),1\leq
j\leq
n,j\not=i\}$.
%%%%%%%%%%%%%%%%
The empirical procedure of \cite{YJ} could be used to get the
alternative data-driven global bandwidth
\[
h_{yj}=h_{\mathrm{cv}} \biggl(\frac{\beta_n(1-\beta_n)}{\phi(\Phi^{-1}(\beta_n))^2}
\biggr)^{1/5},
\]
where $\phi$ and $\Phi$ stand, respectively, for the standard normal
density and distribution functions.
%%%%%%%
However, the use of the `RP' estimators
$\tilde{q}^{\mathrm{RP},i}_n(\beta_n|x)$ and $\hat{\gamma
}^{\mathrm{RP},i}_n(x)$,
for $i=1,2$, requires in addition the selection of an appropriate order
$\alpha_n$.
To simplify the discussion, we set $\alpha_n$ at $k/n^{\star}_{xh}$,
where the
integer $k$ varies between $1$ and $n^{\star}_{xh}-1$, for each $h\in
\mathcal{H}$.
We also consider the value $J=3$ for which $\hat{\gamma
}^{\mathrm{RP},1}_n(x)\equiv\hat{\gamma}^{\mathrm
{RP},2}_n(x):=\hat{\gamma}^{\mathrm{RP}}_n(x;h,k)$
and
$\tilde{q}^{\mathrm{RP},1}_n(\beta_n|x)\equiv\tilde{q}^{\mathrm
{RP},2}_n(\beta_n|x):=\tilde{q}^{\mathrm{RP}}_n(\beta_n|x;h,k)$,
with $r=1/J$.
%%%
An empirical way to decide what values of $(h,k)$ should one
use to compute the estimates in practice could be the automatic {ad
hoc} data driven-rule employed in~\cite{DFS}.
The main idea is to evaluate first the estimates, for each $x$ in a
chosen grid of values, and then to select the parameters where the
variation of the results is the smallest. This can be achieved in two ways:

\textit{Selecting $h$ and $k$ separately.}

\textit{Step} 1. Select a data-driven global bandwidth $h$, for example,
$h_{\mathrm{cv}}$ or $h_{yj}$.

\textit{Step} 2. Evaluate
$\tilde{q}^{\mathrm{RP}}_n(\beta_n|x;h,k)$ at $k=1,\ldots,n^{\star
}_{xh}-1$. Then
compute the standard deviation of the estimates over a \textit{`window'}
of (say
$[\sqrt{n^{\star}_{xh}}]$) successive values of $k$. The value of $k$ where
this standard deviation is minimal defines the desired parameter.

%In other words, even though any data-driven method has to select $h$
%and $k$
%simultaneously, we want to test here the performance of the estimate
%$\tilde{q}^{RP}_n(\beta_n|x;h_{\mathrm{cv}},k_{xh_{\mathrm{cv}}})$.
The same considerations evidently apply to $\hat{\gamma}^{\mathrm
{RP}}_n(x)$ and
to the `benchmark' estimators
$\hat{\gamma}^{\mathrm{GP}}_n(x):=\hat{\gamma}^{\mathrm
{GP}}_n(x;h,k)$ and
$\hat{q}^{\mathrm{GP}}_n(\beta_n|x):=\hat{q}^{\mathrm{GP}}_n(\beta_n|x;h,k)$, defined in
Section~\ref{BG}, with the covariate dependent threshold being
$u_x:=Y_{(n^{\star}_{xh}-k)}^x$, and $Y_{(1)}^x\leq\cdots\leq
Y_{(n^{\star}_{xh})}^x$ being the sequence of ascending order
statistics corresponding to the $Y_i$'s such that
$|X_i-x|\leq h$.

The main difficulty when employing such a separate choice of $h$ and
$k$ is that both
$\tilde{q}^{\mathrm{RP}}_n(\beta_n|x;h_{\mathrm{cv}},k)$ and
$\hat{q}^{\mathrm{GP}}_n(\beta_n|x;h_{\mathrm{cv}},k)$, respectively,
$\hat{\gamma}^{\mathrm{RP}}_n(x;h_{\mathrm{cv}},k)$ and $\hat{\gamma
}^{\mathrm{GP}}_n(x;h_{\mathrm{cv}},k)$, as
functions of $k$ may be so unstable that reasonable values of $k$
(which would
correspond to the true value of $q(\beta_n|x)$, respectively, $\gamma
(x)$) may be hidden in the graphs.
In result, the estimators may exhibit considerable volatility as
functions of
$x$ itself.

%f2 #&#
\begin{figure}

\includegraphics{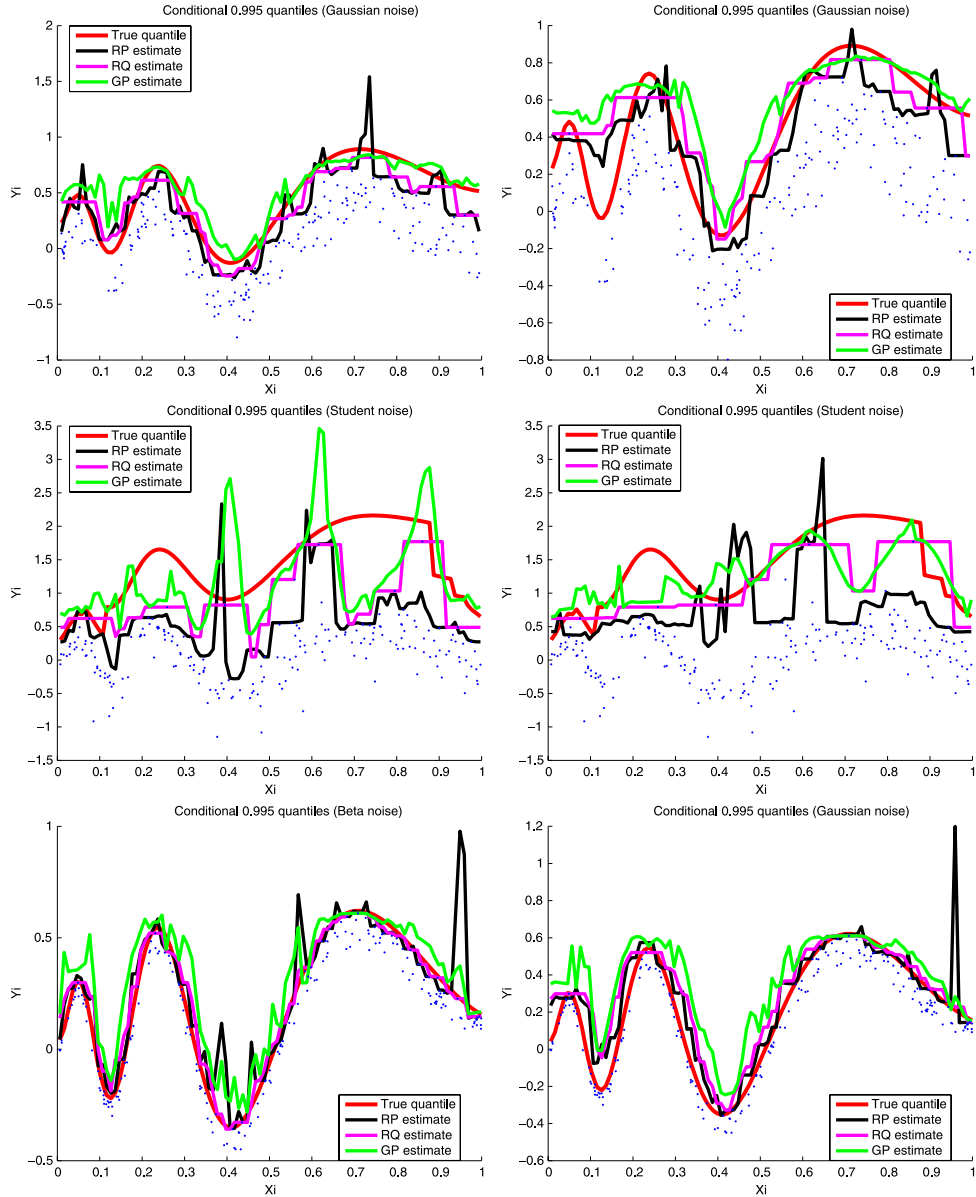}

\caption{Separate parameters' selection using the same simulated
samples as in Figure~\protect\ref{figA}.
From left to right, the used bandwidth is $h_{\mathrm{cv}}$, $h_{yj}$.
From top to bottom, $Y|X$ is Gaussian, Student, Beta.
The true quantile function $q(\beta_n|\cdot)$ in red with
$\beta_n=1/n$. Its estimators $\hat{q}_n(\beta_n|\cdot)$ in magenta,
$\tilde{q}^{\mathrm{RP},1}_n(\beta_n|\cdot)\equiv\tilde
{q}^{\mathrm{RP},2}_n(\beta_n|\cdot)$
in black with $r=1/J$ and $J=3$, and $\hat{q}^{\mathrm{GP}}_n(\beta_n|\cdot)$ in
green. The observations $(X_i,Y_i)$ are depicted as blue points.}
\label{figdatadrivenruleseparate}
\end{figure}

A typical realization is shown in
Figure~\ref{figdatadrivenruleseparate} when the bandwidths
$h_{\mathrm{cv}}$ (left
panels) and $h_{yj}$ (right panels) are used in Step~1. It may be seen that
the method affords reasonable estimates in both Gaussian and Beta error
models regarding the
difficult curvature of the extreme quantile regression and the very small
sample size $n=200$.
%, and the arbitrary choice $(3,1/3)$ of the parameters $(J,r)$.
However, it seems that the method fails in the case of Student noise, where
the superiority of $\tilde{q}^{\mathrm{RP}}_n(\beta_n|\cdot)$ over both
$\hat{q}^{\mathrm{GP}}_n(\beta_n|\cdot)$ and $\hat{q}_n(\beta_n|\cdot)$,
demonstrated via the Monte Carlo study, is clearly sacrificed.
This failure is probably due to the arbitrary choice $(3,1/3)$ of the
parameters $(J,r)$ in $\tilde{q}^{\mathrm{RP}}_n(\beta_n|\cdot)$.
%Thereby, our \textit{ad hoc} technique might be viewed as an
%exploratory
%tool in %the heavy-tailed case, rather than as a method for final
%analysis.
It might also be seen that, apart from the student error model, the three
estimators $\tilde{q}^{\mathrm{RP}}_n(\beta_n|\cdot)$, $\hat
{q}_n(\beta_n|\cdot
)$ and
$\hat{q}^{\mathrm{GP}}_n(\beta_n|\cdot)$ point toward similar results.

\textit{Selecting $h$ and $k$ simultaneously.}

\textit{Step} 1. For each $h\in\mathcal{H}$, proceed to Step~2 described
in the
separate parameters' selection.
Set the value of $k$ where the standard deviation is minimal to be $k_{xh}$
and calculate the corresponding estimate
$\tilde{q}^{\mathrm{RP}}_n(\beta_n|x;h,k_{xh})$.

\textit{Step} 2. Compute the standard deviation of the estimates
$\tilde{q}^{\mathrm{RP}}_n(\beta_n|x;h,k_{xh})$ over a window of
(say $10$)
successive values of $h$. Select the bandwidth where the standard
deviation is minimal and then evaluate the corresponding estimate.

%%%%%%%%%%%%%%%%%%%
%
%f3 #&#
\begin{figure}

\includegraphics{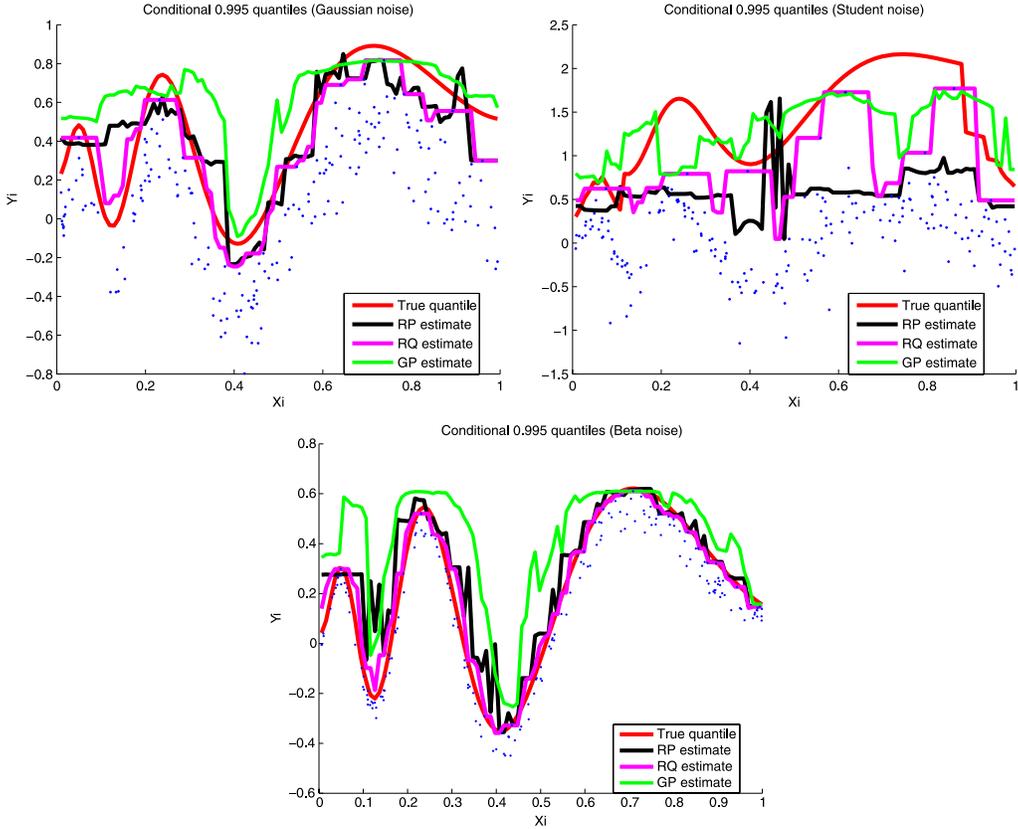}

\caption{Simultaneous parameters' selection using the same simulated
samples as above.
From left to right and from top to bottom, $Y|X$ is Gaussian, Student, Beta.
The true quantile function $q(\beta_n|\cdot)$ in red with
$\beta_n=1/n$. Its estimators $\hat{q}_n(\beta_n|\cdot)=\hat
{q}_n(\beta_n|\cdot,h_{\mathrm{cv}})$ in magenta,
$\tilde{q}^{\mathrm{RP},1}_n(\beta_n|\cdot)\equiv\tilde
{q}^{\mathrm{RP},2}_n(\beta_n|\cdot)$
in black with $(r,J)=(1/3,3)$, and $\hat{q}^{\mathrm{GP}}_n(\beta_n|\cdot)$ in
green. The $(X_i,Y_i)$'s are depicted as blue points.}
\label{figdatadrivenrulesimultaneous}
\end{figure}

In our simulations, we used a refined grid $\mathcal{H}$ of $50$ points
between $\min(h_{\mathrm{cv}},h_{yj}-h_{\mathrm{cv}})$ and $h_{yj}+2h_{\mathrm{cv}}$. Any other limit
bounds of $\mathcal{H}$ could of course be chosen near $h_{\mathrm{cv}}$ below
and near
$h_{yj}$ above.
See Figure~\ref{figdatadrivenrulesimultaneous} for a typical realization
in each simulated scenario.
Here also the method is not without disadvantage as can be seen from
the case
of Student noise, where good results require a large sample size.
%The `RQ' estimator
%$\hat{q}_n(\beta_n|x)$ has, however, a nice
%behaviour especially in the Gaussian and Beta error models. The `GP'
%estimator $\hat{q}^{\mathrm{GP}}_n(\beta_n|x)$ seems to be the winner
%in the
%Student er%ror model.
%%%%%%%%%%%%%%%%%%%
%

%as
%in
% Figure~\ref{figA}.
%Left (right) panels: simultaneous (separate) parameters' selection.
%From top to bottom, $Y|X$ is Gaussian, Student, Beta.
%The true quantile function $q(\beta_n|\cdot)$ in red with
% $\beta_n=1/n$. Its estimators $\hat{q}_n(\beta_n|\cdot)$ in magenta,
% $\tilde{q}^{\mathrm{RP},1}_n(\beta_n|\cdot)\equiv\tilde{q}^{
% in black with $r=1/J$ and $J=3$, and $\hat{q}^{\mathrm{GP}}_n(\beta_n|
%in
% green. The observations $(X_i,Y_i)$ are depicted as blue points.}

%s4.4 #&#
\subsection{Concluding remarks}
\textit{Monte Carlo evidence}.
The experiments indicate that $\hat{q}_n(\beta_n|x)$
is efficient relative to the modern smoothing estimator $\hat
{q}^{\mathrm{GP}}_n(\beta_n|x)$ first introduced in \cite{BEG,BES}.
The simulations also indicate that the performance of the alternative
estimator $\tilde{q}_n^{\mathrm{RP},1}(\beta_n|x)$ is quite
remarkable in comparison with its analog
$\tilde{q}_n^{\mathrm{RP},2}(\beta_n|x)$, at least in terms of MSE.
In comparison with $\hat{q}_n(\beta_n|x)$ and $\hat{q}^{\mathrm
{GP}}_n(\beta_n|x)$, the
variability and the bias of both $\tilde{q}_n^{\mathrm{RP},1}(\beta_n|x)$ and $\tilde{q}_n^{\mathrm{RP},2}(\beta_n|x)$ are quite
respectable and it seems that the heavier is the conditional tail, the better
the estimators $\tilde{q}_n^{\mathrm{RP}}(\beta_n|x)$ are.
It should be also clear that $\tilde{q}_n^{\mathrm{RP},1}(\beta_n|x)$ and $\tilde{q}_n^{\mathrm{RP},2}(\beta_n|x)$ can be improved
appreciably by tuning the choice of the parameters $J$ and $r$.

\textit{Tuning parameters selection in practice.}
The simulations worked reasonably well for our `ad hoc' selection
methods except for the heavy-tailed case, corresponding
to generally severe events.
A sensible practice would be to verify whether the resulting `RP', `GP' and
`RQ' estimators point toward similar conclusions: the hard question of
how to
pick out the smoothing parameters $(h,k)$ simultaneously in an optimal way
might thus become less urgent. In contrast, if the estimators look
clearly different, this might diagnose a heavy-tailed conditional
distribution with a great variability in severity:
thereby our technique might be viewed as an exploratory tool, rather
than as a method for final analysis.
%Nevertheless, in the latter case, our
%simultaneous-selection method seems to provide respectable `GP'
%estimates.
Doubtlessly, further work to define a concept of selecting appropriate
values for the crucial parameters $(J,r)$ in the $\tilde{q}_n^{\mathrm
{RP}}(\beta_n|x)$ estimators will
yield new refinements.

\textit{The case of multiple covariates.}
We have discussed the asymptotic distributional properties of both estimators
$\hat{q}_n(\beta_n|x)$ and $\tilde{q}_n^{\mathrm{RP}}(\beta_n|x)$
in detail in Sections \ref
{notation} and \ref{results} for multiple regressors $X\in\mathbb{R}^p$,
but our contributions are
probably only of a theoretical value in the case $p>1$.
Indeed, as in the ordinary setting where the quantile order does not
depend on
the sample size, the kernel-smoothing method suffers from the `curse of
dimensionality'. In our setting of extreme quantile regression, the
curse is
exacerbated by several degrees of magnitude and drastically increases in
higher dimensions. To overcome this vexing defect, one can use dimension
reduction techniques such as ADE (Average Derivative Estimator), see for
instance~\cite{HARDLE}. Nevertheless, the theoretical properties of such
methods are not yet established in the extreme-value framework.

%As may be seen from the figure, this \textit{ad hoc} prescription
%results
%in
%satisfactory estimates in terms of both bias and stability, although
%they do
%not enjoy the desired optimality.

%s5 #&#
\section{Data example}\label{data}
%%%%%%%%%%%%%%%%%%%%%%%%%%%%%%%%%%%%%%%%%%%%%%%%%%%%%%%%%%%%%%%%%%%%%%%%%%%%%%%%

Data on $123$ American electric utility companies were collected and
the aim
is to investigate the economic efficiency of these companies (see,
e.g., \cite{Gre}). A possible way to measure this efficiency is by looking
at the maximum level of produced goods which is attainable for a given level
of inputs-usage. From a statistical point of view, this problem translates
into studying the upper boundary of the set of possible inputs $X$ and
outputs $Y$, the so-called cost/econometric frontier in production theory.
Hendricks and Koenker~\cite{HEN} stated: ``In the econometric
literature on
the estimation of production technologies, there has been
considerable\vadjust{\goodbreak}
interest in estimating so called frontier production models that correspond
closely to models for extreme quantiles of a stochastic production surface''.
The present paper may be viewed as the first `purely' nonparametric
work to actually investigate theoretically the idea of Hendricks and Koenker.

%
%f4 #&#
\begin{figure}

\includegraphics{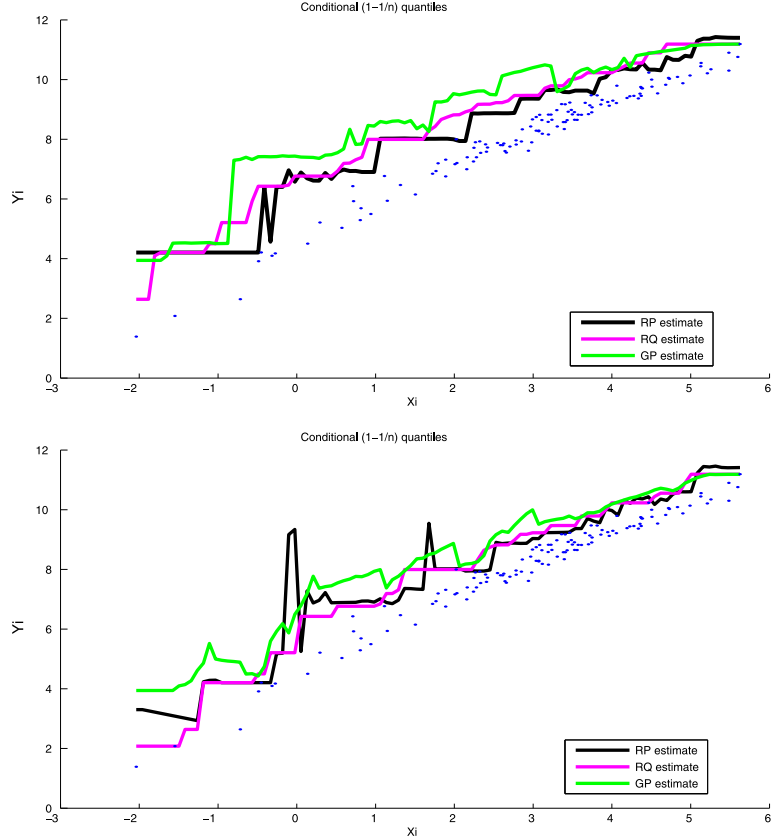}

\caption{Scatterplot of the American electric utility
data (blue points). Results obtained via the simultaneous (top) and
separate (bottom) selection methods for $\beta_n=1/n$: the estimator
$\hat{q}_n(\beta_n|\cdot)$ is drawn in magenta with $h_{yj}$ (top) and
$h_{\mathrm{cv}}$ (bottom), $\hat{q}^{\mathrm{GP}}_n(\beta_n|\cdot)$ in
green and
$\tilde{q}^{\mathrm{RP},1}_n(\beta_n|\cdot)$ in black with $(r,J)=(1/3,3)$.}
\label{greene1}
\end{figure}

For our illustration purposes, we used the measurements on
the variable $Y=\log(Q)$, with $Q$ being the production output of a
firm, and
the variable $X=\log(C)$, with $C$ being the total cost involved in the
production.
Figure~\ref{greene1} shows the $n=123$ observations,
together with estimated extreme conditional quantiles
$\hat{q}_n(\beta_n|x)$,
$\hat{q}^{\mathrm{GP}}_n(\beta_n|x)$
and $\tilde{q}_n^{\mathrm{RP},1}(\beta_n|x)=\tilde{q}_n^{\mathrm
{RP},2}(\beta_n|x)$ at $r=1/J$ with $J=3$.
Given the small sample size, it was enough to use
$\beta_n=1/n$
in describing the conditional distribution tails.
For selecting the window width $h$ and the number $k$ of extremes, we
maintained the automatic empirical data-driven rules described above.
%%%%%%%%%%%%%%%%%%%
%However, it may be noticed that the use of $h_{\mathrm{min}}:=\max_{1\leq
% i<n}\left(X_{(i+1)}-X_{(i)}\right)$ as a lower bound of $\mathcal{H}$
%for
%this particular dataset will result in large values of $h$ due to the
%design
%sparseness in the left tail. For instance, the bandwidth
%$h_{\mathrm{cv}}$ computed via the cross-validation criterion was found to be
%identical
%to $h_{\mathrm{min}}\equiv0.8305$. To avoid misleading estimates we used the
%more reas%onable lower bound
%$h_{\mathrm{min}}:=\min_{1\leq i<n}\left(X_{(i+1)}-X_{(i)}\right)\equiv
%1.0571e-04$,
%which results in $h_{\mathrm{cv}}=0.5667$ and the graphs in Figure~
%%%%%%%%%%%%%%%%%%%
It appears that the extreme conditional quantile estimates are similar
for both simultaneous (top) and separate (bottom) selection methods.
%%%%%%%%%%%%%%%%%%%
Following their evolution,
it may be seen that the American electric utility data do not
correspond to
the situation hoped for by the practitioners of a heavily short-tailed
production process. Indeed, one may distinguish between two
different behaviors of the extreme regression quantiles: They indicate
a short-tailed conditional distribution for companies working at (transformed)
input-factors larger than, say, $2$. In contrast, the tail distribution for
the smallest companies with inputs $X_i< 2$ seems to be moderately heavy.
%%%%%%%%%%%%%%%%%%%
%graphed in Figure~\ref{greene2},
%which was found to vary modestly with the sample boundary.
%As is to be expected from the Monte Carlo evdience, the `benchmark'
%estimator %$\hat{\gamma}^{\mathrm{GP}}_n(x)$
%has an admirable and more stable
%behaviour than $\gamrpicku=\gamrpickd$ for the chosen parameters
%$(r,J)$.
%Besides its large variability, the latter estimator suffers
%drastically from % the design sparseness in the left tail.
% }
%%%%%%%%%%%%%%%%%%%
%As a matter of fact, good extreme-value index estimates
%in this context of production data should vary between $-1$ and $0$
%with the
%following intuitive interpretation~\cite{DFS}:
% When $\gamma(x)<-\{p+1\}^{-1}$, the conditional density rises up to
%infinity
% as it
%%%%%%%%%%%%%%%%%%%%
%%%%%%%%%%%%%%%%%%%%
%which would correspond to an ideal production activity.
%The situation hoped for by the practitioners is, at most, the case
%$\gamma(x)=-\{p+1\}^{-1}$ where the density is strictly positive at
%the fronti%er.
%The density of data often decays to zero smoothly as
%it
%%%%%%%%%%%%%%%%%%%%
%%%%%%%%%%%%%%%%%%%
% which corresponds to the case $\gamma(x)>-\{p+1\}^{-1}$.
%%%%%%%%%%%%%%%%%%%
%$\hat{\gamma}^{\mathrm{GP}}_n(x)\geq-0.5$,}
%%%%%%%%%%%%%%%%%%%%
%and so the American electric utility
%%%%%%%%%%%%%%%%%%%%%%%%%%%%%%%%%%%%%%%%
Therefore, the theoretical economic assumption that producers should
operate on the upper boundary of the joint support of $(X,Y)$ rather
than on
its interior is clearly not fulfilled here, revealing a certain lack of
production performance in this sector of activity.
%%%%%%%%%%%%%%%%%%%%%%%%%%%%%%%%%%%%%%%%%
%It may be also noticed that, due to the small sample size here, the
%estimator
%$\gamrpicku$ exhibits some volatility (for the smallest
%inputs) exceeding the upper
%limit bound $0$, but it remains quite reasonable.
%%%%%%%%%%%%%%%%%%%
The estimated graph of $\tilde{q}_n^{\mathrm{RP},1}(\beta_n|x)$,
$\hat{q}_n(\beta_n|x)$ or
$\hat{q}^{\mathrm{GP}}_n(\beta_n|x)$ might be interpreted as the set
of the most
efficient firms.
%%%%%%%%%%%%%%%%%%%
It is then clear
that the firms achieve significantly lesser output than that predicted
by the
extremal quantile frontiers. This indicates a relative economic inefficiency
especially in the population of the (sparse) smallest companies in
terms of inputs.\looseness=-1\vspace*{-3pt}
%%%%%%%%%%%%%
%As shown by Gijbels \textit{et al.} (1999), the parametric translog and
%Cobb-Doug%las models, commonly used in applied
%econometrics and which describe a linear model and a quadratic model
%of the
%output $Y$ in terms of the input $X$, are not so appropriate for the
%American %electric utility data.
%%%%%%%%%%%%%%%%%%%
%%%%%%%%%%%%%%%%%%%

%%%%%%%%%%%%%%%

%%%%%%%%%%%%%%%%%%%%
% \begin{center}
% \end{center}
% \vspace*{-0.4 cm}
%
%simultaneous (l-h.%s) and separate (r-h.s) selection methods. The
%inputs $X_i$ are depicted as
% blue points, the estimator $\gamrpicku=\gamrpickd$ is drawn in black
%and
% $\hat{\gamma}^{\mathrm{GP}}_n(x)$ in green.}
% \label{greene2}
%%%%%%%%%%%%%%%

\begin{appendix}\label{appendix}
%%%%-----------------------------------------------------------------------------
%s6 #&#
\section*{Appendix: Proofs}
%%%%--------------------------------------------------------------------------

%s6.1 #&#
\subsection{Preliminary results}
%%%%--------------------------------------------------------------------------

We begin with a homogeneous property of the quantile function.
%
%le1 #&#
\begin{Lem}
\label{lemrapquant}
Suppose \textup{(A.1)} holds. If $\alpha_n\to0$ as $n\to\infty$, then,
\[
\lim_{n\to\infty} \frac{q(\xi\alpha_n|x)}{q(\alpha_n|x)} = \xi^{-(\gamma(x) \vee0)}
\]
for all $\xi>0$.
\end{Lem}
\begin{pf} From (\ref{F1quant}), we have
\[
\frac{q(\xi\alpha_n|x)}{q(\alpha_n|x)} = 1 + K_{\gamma(x)}(1/\xi) \frac
{a(q(\alpha_n|x)|x)}{q(\alpha_n|x)} \bigl(1+
\mathrm{o}(1) \bigr)
\]
and the conclusion follows using~(\ref{clef}).
\end{pf}

%%%%-------------------------------------------------------------------------

The following lemma states that the convergence in~(\ref{Fu})
is uniform.

%le2 #&#
\begin{Lem}
\label{lemunif}
Under \textup{(A.1)}, if $z_n(x) \uparrow y_F(x)$ as $n \to\infty$, then for all
sequence of functions $t_n(x)$ such that
$t_n(x) \to t_0(x)$ as $n \to\infty$ where $t_0(x)$ is such that there
exists $\eta>0$ for which $1+\gamma(x)t_0(x) \geq\eta$ then,
\[
\lim_{n \to\infty} \frac{\baF(z_n(x)+t_n(x)a(z_n(x)|x)|x)}{\baF
(z_n(x)|x)} = \frac{1}{K_{\gamma(x)}^{-1}(t_0(x))}.
\]
\end{Lem}

\begin{pf} Since $t_n(x) \to t_0(x)$ as $n \to\infty$, for all
$\varepsilon_1>0$ such that $|\gamma(x)|\varepsilon_1<\eta$, there
exists $N_1\geq0$ such that\vadjust{\goodbreak} for all $n \geq N_1$, $t_0(x)-\varepsilon_1 \leq t_n(x) \leq t_0(x)+\varepsilon_1$. Since $a(z_n(x)|x)>0$ and
$\baF(\cdot|x)$ is a decreasing function, we have:
\begin{eqnarray*}
\frac{\baF(z_n(x)+(t_0(x)+\varepsilon_1)a(z_n(x)|x)|x)}{\baF
(z_n(x)|x)} &\leq& \frac{\baF(z_n(x)+t_n(x)a(z_n(x)|x)|x)}{\baF(z_n(x)|x)}
\\
&\leq& \frac{\baF(z_n(x)+(t_0(x)-\varepsilon_1)a(z_n(x)|x)|x)}{\baF
(z_n(x)|x)}.
\end{eqnarray*}
Now, since $|\gamma(x)|\varepsilon_1<\eta$, it is easy to check that
$1+\gamma(x)(t_0(x)+\varepsilon_1)\wedge1+\gamma
(x)(t_0(x)-\varepsilon_1) >0$. Hence, under~(A.1), for all
$\varepsilon_2 >0$, there exists
$N_2\geq0$ such that for all $n \geq N_2$
\[
\frac{1-\varepsilon_2}{K_{\gamma(x)}^{-1}(t_0(x)+\varepsilon_1)}\leq\frac{\baF(z_n(x)+t_n(x)a(z_n(x)|x)|x)}{\baF(z_n(x)|x)} \leq\frac{1+\varepsilon_2}{K_{\gamma(x)}^{-1}(t_0(x)-\varepsilon_1)}.
\]
Since $\varepsilon_1$ and $\varepsilon_2$ can be chosen arbitrarily
small, this concludes the proof.
\end{pf}

%%%%--------------------------------------------------------------------------
Let us remark that the kernel estimator~(\ref{defestproba})
can be rewritten as $\habaF(y|x)=\hat\psi_n(y,x)/\hat g_n(x)$
where
\[
\hat\psi_n(y,x)= \frac{1}{n} \sum
_{i=1}^n K_h(x-X_i) \I
\{Y_i>y\}
\]
is an estimator of $\psi(y,x)= \baF(y|x)g(x)$ and
$\hat g_n(x)$ is the classical kernel estimator of the p.d.f. $g(x)$
defined by:
\[
\hat g_n(x) = \frac{1}{n} \sum_{i=1}^n
K_h(x-X_i).
\]
Lemma~\ref{lemdensite} gives standard results on the
kernel estimator
%%%%%%%%%%%%%%%%%%%
(see \cite{PARZEN} for a proof)
%%%%%%%%%%%%%%%%%%
whereas Lemma~\ref{lempsi} is dedicated to the asymptotic
properties of $\hat\psi_n(y,x)$.
%%%%--------------------------------------------------------------------------

%le3 #&#
\begin{Lem}
\label{lemdensite}
Suppose \textup{(A.2)}, \textup{(A.3)} hold. If $nh^p\to\infty$, then, for all $x\in\R^p$
such that \mbox{$g(x)>0$},
\begin{longlist}[(ii)]
\item[(i)] $\E(\hat g_n(x)-g(x)) = \mathrm{O}(h)$,
\item[(ii)] $\operatorname{var}(\hat g_n(x))= \frac{g(x) \|K\|_2^2}{nh^p}(1+\mathrm{o}(1))$.
\end{longlist}
\end{Lem}

Therefore, under the assumptions of the above lemma, $\hat g_n(x)$
converges to $g(x)$ in probability.

%%%%--------------------------------------------------------------------------
Let us introduce some further notation. In the following,
$y_n(x)$ is a sequence such that $y_n(x) \uparrow y_F(x)$
and $y_{n,j}(x) = y_n(x) + K_{\gamma(x)}(1/\tau_j) a(y_n|x)(1+\mathrm{o}(1))$
for all $j=1,\ldots,K$.
Recall that $0<\tau_J< \cdots< \tau_2 < \tau_1 \leq1$.
Moreover, the oscillations of the c.s.f. are controlled by
\[
\omega_n(x) := \max_{j=1,\ldots,J}\sup_{x'\in B(x,h)} \biggl\llvert
\frac
{\baF
(y_{n,j}(x)|x')}{\baF(y_{n,j}(x)|x)}-1 \biggr\rrvert.
\]
%
%le4 #&#
\begin{Lem}
\label{lempsi}
Suppose \textup{(A.1)}, \textup{(A.2)} and \textup{(A.3)} hold and let $x\in\R^p$ such that $g(x)>0$.
If $\omega_n(x) \to0$ and $nh^p\baF(y_n(x)|x) \to\infty$ then,
\begin{longlist}[(ii)]
\item[(i)] $\E(\hat\psi_n(y_{n,j}(x),x)) = \psi(y_{n,j}(x),x)
(1+\mathrm{O}(\omega_n(x))+ \mathrm{O}(h) )$,
for $j=1,\ldots, J$.
\item[(ii)] The random vector
\[
\biggl\{ \sqrt{nh^p\psi \bigl(y_n(x),x \bigr)} \biggl(
\frac{\hat\psi_n(y_{n,j}(x),x)
-\E(\hat\psi_n(y_{n,j}(x),x))}{\psi(y_{n,j}(x),x)} \biggr) \biggr\}_{j=1,\ldots,J}
\]
is asymptotically Gaussian, centered, with covariance matrix
$\|K\|_2^2 V$ where $V_{j,j'} = \tau^{-1}_{j\wedge j'}$
for $(j,j')\in\{1,\ldots, J\}^2$.
\end{longlist}
\end{Lem}

\begin{pf} (i) Since the $(X_i,Y_i)$, $i=1,\ldots,n$ are identically
distributed, we have
\begin{eqnarray*}
\E \bigl(\hat\psi_n \bigl(y_{n,j}(x),x \bigr) \bigr)&=&
\int_{\R^p} K_h(x-t) \baF \bigl(y_{n,j}(x)|t
\bigr) g(t) \,\mathrm{d}t\\
 &=& \int_{S} K(u) \baF
\bigl(y_{n,j}(x)|x-hu \bigr) g(x-hu) \,\mathrm{d}u,
\end{eqnarray*}
under (A.3). Let us now consider
\renewcommand{\theequation}{\arabic{equation}}
%e8 #&#
%e9 #&#
\begin{eqnarray}\label{eq1}
&&\bigl|\E \bigl(\hat\psi_n \bigl(y_{n,j}(x),x \bigr) \bigr) -
\psi \bigl(y_{n,j}(x),x \bigr)\bigr|
\nonumber
\\[-8pt]
\\[-8pt]
\nonumber
&&\quad  \leq \baF
\bigl(y_{n,j}(x)|x \bigr) \int_S K(u)
\bigl|g(x-hu)-g(x)\bigr|\,\mathrm{d}u
\\
\label{eq2}
& &\qquad{} + \baF \bigl(y_{n,j}(x)|x \bigr) \int_S
K(u) \biggl\llvert\frac
{\baF(y_{n,j}(x)|x-hu)}{\baF(y_{n,j}(x)|x)} - 1 \biggr\rrvert g(x - hu) \,\mathrm{d}u.
\end{eqnarray}
Under (A.2), and since $g(x)>0$, we have
%
%e10 #&#
\begin{equation}
\label{eq4} (\ref{eq1}) \leq\baF \bigl(y_{n,j}(x)|x \bigr)
c_g h \int_S \,\mathrm{d}(u,0) K(u) \,
\mathrm{d}u = \psi \bigl(y_{n,j}(x),x \bigr) \mathrm{O}(h),
\end{equation}
and, in view of (\ref{eq4}),
%
%e11 #&#
\begin{eqnarray}\label{eq3}
(\ref{eq2}) &\leq& \baF \bigl(y_{n,j}(x)|x \bigr)
\omega_n(x) \int_S K(u) g(x-hu) \,\mathrm{d}u =
\baF \bigl(y_{n,j}(x)|x \bigr) \omega_n(x) g(x) \bigl(1+
\mathrm{o}(1) \bigr)
\nonumber
\\[-8pt]
\\[-8pt]
\nonumber
&=& \psi \bigl(y_{n,j}(x) ,x \bigr) \omega_n(x) \bigl(1+
\mathrm{o}(1) \bigr) .
\end{eqnarray}
Combining (\ref{eq4}) and (\ref{eq3}) concludes the first part of the
proof.

(ii) Let $\beta\neq0$ in $\R^J$,
$\Lambda_n(x)= (nh^p\psi(y_n(x),x))^{-1/2}$, and consider the random
variable
\begin{eqnarray*}
\Psi_n&=&\sum_{j=1}^J
\beta_j \biggl(\frac{\hat\psi_n(y_{n,j}(x),x)
-\E(\hat\psi_n(y_{n,j}(x),x))}{\Lambda_n(x)\psi
(y_{n,j}(x),x)} \biggr)
\\
&=&\sum_{i=1}^n \frac{1}{n \Lambda_n(x)} \Biggl
\{ \sum_{j=1}^J \frac
{\beta_j K_h(x-X_i) \I\{Y_i\geq y_{n,j}(x)\} }{\psi(y_{n,j}(x),x)}\\
&&\hspace*{46pt}\quad{} - \E
\Biggl( \sum_{j=1}^J \frac{\beta_j K_h(x-X_i)\I\{Y_i\geq
y_{n,j}(x)\}
}{\psi(y_{n,j}(x),x)}
\Biggr) \Biggr\}
\\
&:=& \sum_{i=1}^n Z_{i,n}.
\end{eqnarray*}
Clearly, $\{Z_{i,n}, i=1,\ldots,n\}$ is a set of centered, independent
and identically distributed random variables with variance
\[
\operatorname{var}(Z_{i,n}) = \frac{1}{n^2 h^{2p}\Lambda_n^2(x)} \operatorname{var} \Biggl( \sum
_{j=1}^J \beta_j K \biggl(
\frac{x-X_i}{h} \biggr) \frac{\I\{Y_i\geq y_{n,j}(x)\} }{\psi
(y_{n,j}(x),x)} \Biggr) = \frac{1}{n^2 h^p\Lambda_n^2(x)}
\beta^t B \beta,
\]
where $B$ is the $J\times J$ covariance matrix
with coefficients defined for $(j,j')\in\{1,\ldots,J\}^2$ by
\begin{eqnarray*}
B_{j,j'}&=& \frac{A_{j,j'}}{\psi(y_{n,j}(x),x)\psi(y_{n,j'}(x),x)},
\\
A_{j,j'} &=& \frac{1}{h^p}\operatorname{cov} \biggl( K \biggl(
\frac
{x-X}{h} \biggr) \I \bigl\{Y\geq y_{n,j}(x) \bigr\} , K
\biggl( \frac{x-X}{h} \biggr) \I \bigl\{Y\geq y_{n,j'}(x) \bigr\}
\biggr)
\\
&=& \|K\|_2^2 \E \biggl(\frac{1}{h^p} Q \biggl(
\frac{x-X}{h} \biggr) \I \bigl\{Y\geq y_{n,j}(x) \vee
y_{n,j'}(x) \bigr\} \biggr)
\\
&&{}- h^p \E \bigl( K_h(x-X) \I \bigl\{Y\geq
y_{n,j}(x) \bigr\} \bigr) \E \bigl( K_h(x-X) \I \bigl\{ Y
\geq y_{n,j'}(x) \bigr\} \bigr),
\end{eqnarray*}
with $Q(\cdot):=K^2(\cdot)/\|K\|^2_2$ also satisfying assumption (A.3).
As a consequence, the three above expectations are of the same nature.
Thus, remarking that, for $n$ large enough,
$y_{n,j}(x) \vee y_{n,j'}(x)=y_{n,j\vee j'}(x)$,
part (i) of the proof implies
\begin{eqnarray*}
A_{j,j'}& =& \|K\|^2_2 \psi
\bigl(y_{n,j\vee j'}(x),x \bigr) \bigl(1+\mathrm{O} \bigl(\omega_n(x)
\bigr)+\mathrm{O}(h) \bigr)
\\
&&{}- h^p \psi \bigl(y_{n,j}(x),x \bigr)\psi
\bigl(y_{n,j'}(x),x \bigr) \bigl(1+\mathrm{O} \bigl(\omega_n(x)
\bigr)+\mathrm{O}(h) \bigr)
\end{eqnarray*}
leading to
\begin{eqnarray*}
B_{j,j'} &= &\frac{\|K\|^2_2}{\psi(y_{n,j\wedge j'}(x),x)} \bigl(1+\mathrm{O} \bigl(
\omega_n(x) \bigr)+\mathrm{O}(h) \bigr) - h^p \bigl(1+
\mathrm{O} \bigl(\omega_n(x) \bigr)+\mathrm{O}(h) \bigr)
\\
&=& \frac{\|K\|^2_2}{\psi(y_{n,j\wedge j'}(x),x)} \bigl(1+\mathrm {o}(1) \bigr),
\end{eqnarray*}
since $\psi(y_{n,j\wedge j'}(x),x)\to0$ as $n\to\infty$.
Now, from Lemma~\ref{lemunif},
\[
\lim_{n\to\infty} \frac{\psi(y_{n,j\wedge j'}(x),x)}{\psi
(y_{n}(x),x)} = \frac{1}{K_{\gamma(x)}^{-1}(K_{\gamma(x)}(1/\tau_{j
\wedge j'}))}=
\tau_{j \wedge j'}
\]
entailing
\[
B_{j,j'} = \frac{\|K\|^2_2 V_{j,j'}}{\psi(y_n(x),x)} \bigl(1+\mathrm {o}(1) \bigr),
\]
and therefore,
$
\operatorname{var}(Z_{i,n}) \sim\|K\|^2_2 \beta^t V \beta/n,
$
for all $i=1,\ldots,n$. As a preliminary conclusion, the
variance of $\Psi_n$ converges to $\|K\|^2_2 \beta^t V \beta$.
Consequently, Lyapounov criteria for the asymptotic normality of sums of
triangular arrays reduces to
$
\sum_{i=1}^{n}\E\llvert Z_{i,n}\rrvert^{3}
=n\E\llvert Z_{1,n}\rrvert^{3}\rightarrow0.
$
Remark that $Z_{1,n}$ is a bounded random variable:
\begin{eqnarray*}
|Z_{1,n}| &\leq&\frac{ 2 \|K\|_\infty\sum_{j=1}^J |\beta_j|} {
n \Lambda_n(x) h^p \psi(y_{n,J},x)}\\
&=& \frac{2}{\tau_J}\|K
\|_\infty\sum_{j=1}^J |
\beta_j|\Lambda_n(x) \bigl(1+\mathrm{o}(1) \bigr)
\end{eqnarray*}
and thus,
\begin{eqnarray*}
n\E\llvert Z_{1,n}\rrvert^{3} &\leq& \frac{2}{\tau_J} \|K
\|_\infty\sum_{j=1}^J n|
\beta_j| \Lambda_n(x) {\operatorname{var}}(Z_{1,n})
\bigl(1+\mathrm{o}(1) \bigr)
\\
& =& \frac{2}{\tau_J} \|K\|_\infty\|K\|^2_2
\sum_{j=1}^J |\beta_j|
\beta^t V\beta\Lambda_n(x) \bigl(1+\mathrm{o}(1) \bigr)
\to0
\end{eqnarray*}
as $n\to\infty$. As a conclusion, $\Psi_n$ converges in distribution
to a centered Gaussian random variable with variance $\|K\|^2_2 \beta^t
V \beta$ for all $\beta\neq0$ in $\R^p$. The result is proved.
\end{pf}

Let us now focus on the estimation of small tail probabilities
$\baF(y_n(x)|x)$ when $y_n(x)\uparrow y_F(x)$ as $n\to\infty$.
The following result provides sufficient conditions for
the asymptotic normality of $\habaF(y_n(x)|x)$.

%pr1 #&#
\begin{prop}
\label{thproba}
Suppose \textup{(A.1)}, \textup{(A.2)} and \textup{(A.3)} hold and let $x\in\R^p$ such that
$g(x)>0$. If
$nh^p\baF(y_n(x)|x) \to\infty$ and $nh^{p}\baF(y_n(x)|x) (h \vee
\omega_n(x))^2 \to0$,
then, the random vector
\[
\biggl\{ \sqrt{nh^p\baF \bigl(y_n(x)|x \bigr)} \biggl(
\frac{\habaF(y_{n,j}(x)|x)
}{\baF(y_{n,j}(x)|x)}-1 \biggr) \biggr\}_{j=1,\ldots,J}
\]
is asymptotically Gaussian, centered, with covariance matrix
$\|K\|_2^2/g(x) V$ where $V_{j,j'} = \tau^{-1}_{j\wedge j'}$ for
$(j,j')\in\{1,\ldots, J\}^2$.
\end{prop}

%pa6.1.subsubsection.1 #&#
\begin{pf}
Keeping in mind the notation of Lemma~\ref{lempsi}, the following
expansion holds
%
%e12 #&#
\renewcommand{\theequation}{\arabic{equation}}
\begin{equation}
\label{term123} \Lambda_n^{-1}(x)\sum
_{j=1}^J \beta_j \biggl(
\frac{\habaF
(y_{n,j}(x)|x)}{\baF(y_{n,j}(x)|x)} -1 \biggr) =\frac
{T_{1,n}+T_{2,n}-T_{3,n}}{\hat g_n(x)},
\end{equation}
where
\begin{eqnarray*}
T_{1,n}&=&g(x)\Lambda_n^{-1}(x) \sum
_{j=1}^J \beta_j \biggl(
\frac
{\hat
\psi_n(y_{n,j}(x),x) -
\E(\hat\psi_n(y_{n,j}(x),x)) }{\psi(y_{n,j}(x),x)} \biggr),
\\
T_{2,n}&=&g(x) \Lambda_n^{-1}(x)\sum
_{j=1}^J \beta_j \biggl(
\frac{
\E
(\hat\psi_n(y_{n,j}(x),x)) - \psi(y_{n,j}(x),x) }{\psi
(y_{n,j}(x),x)} \biggr),
\\
T_{3,n}& =& \Biggl( \sum_{j=1}^J
\beta_j \Biggr)\Lambda_n^{-1}(x) \bigl( \hat
g_n(x)-g(x) \bigr).
\end{eqnarray*}
Let us highlight that assumptions $nh^{p}\baF(y_n(x)|x)\omega_n^2(x)
\to0$ and $nh^p \baF(y_n(x)|x) \to\infty$ imply that $\omega_n(x)
\to
0$ as $n\to\infty$.
Thus, from Lemma~\ref{lempsi}(ii), the random term $T_{1,n}$ can be
rewritten as
%
%e13 #&#
\begin{equation}
\label{term1} T_{1,n}= g(x) \|K\|_2 \sqrt{
\beta^t V \beta} \xi_n,
\end{equation}
where $\xi_n$ converges to a standard Gaussian random variable.
The nonrandom term $T_{2,n}$ is controlled with Lemma~\ref{lempsi}(i):
%
%e14 #&#
\begin{equation}
\label{term2} T_{2,n}= \mathrm{O}(\Lambda_n^{-1}(x)
\bigl(h+\Delta \bigl(y_n(x),x \bigr) \bigr)=\mathrm{O} \bigl(
\bigl(nh^p\baF \bigl(y_n(x)|x \bigr) \bigr)^{1/2}
\bigl(h\vee\omega_n(x) \bigr) \bigr)=\mathrm{o}(1).
\end{equation}
Finally, $T_{3,n}$ is a classical term in kernel density estimation,
which can be bounded by Lemma~\ref{lemdensite}:
%
%e15 #&#
\begin{eqnarray}
\label{term3} T_{3,n} & = & \mathrm{O} \bigl(h\Lambda_n^{-1}(x)
\bigr) + \mathrm{O}_P \bigl(\Lambda_n^{-1}(x)
\bigl(nh^p \bigr)^{-1/2} \bigr)
\nonumber
\\[-8pt]
\\[-8pt]
\nonumber
 & = & \mathrm{O} \bigl(n h^{p+2} \baF \bigl(y_n(x)|x \bigr)
\bigr)^{1/2} + \mathrm{O}_P \bigl(\baF \bigl(y_n(x)|x
\bigr) \bigr)^{1/2}=\mathrm{o}_P(1).
\end{eqnarray}
Collecting (\ref{term123})--(\ref{term3}), it follows that
\[
\hat g_n(x) \Lambda_n^{-1}(x)\sum
_{j=1}^J \beta_j \biggl(
\frac
{\habaF
(y_{n,j}(x)|x)}{\baF(y_{n,j}(x)|x)} -1 \biggr) = g(x) \|K\|_2 \sqrt {
\beta^t V \beta} \xi_n +\mathrm{o}_P(1).
\]
Finally, $\hat g_n(x)\toP g(x)$ yields
\begin{eqnarray*}
&&\sqrt{nh^p \baF \bigl(y_n(x)|x \bigr)} \sum
_{j=1}^J \beta_j \biggl(
\frac
{\habaF
(y_{n,j}(x)|x)}{\baF(y_{n,j}(x)|x)} -1 \biggr) \\
&&\quad= \|K\|_2 \sqrt{\frac
{\beta^t V\beta}{g(x)}}
\xi_n +\mathrm{o}_P(1)
\end{eqnarray*}
and the result is proved.
\end{pf}

The last lemma establishes that $K_{\hat\gamma_n(x)}(r_n)$
inherits from the convergence properties of $\hat\gamma_n(x)$.
%
%le5 #&#
\begin{Lem}
\label{lemNAK}
Suppose $\xi_n^{(\gamma)}(x):=\Lambda_n^{-1} (\hat\gamma_n(x)-\gamma
(x))=\mathrm{O}_{\PP}(1)$, where $\Lambda_n\to0$.
Let $r_n\geq1$ or $r_n\leq1$ such that $\Lambda_n\log(r_n)\to0$.
Then,
\[
\Lambda_n^{-1} \biggl( \frac{K_{\hat\gamma_n(x)}(r_n) - K_{\gamma
(x)}(r_n) }{K'_{\gamma(x)}(r_n)} \biggr)=
\xi_n^{(\gamma)}(x) \bigl(1+\mathrm{o}_{\PP}(1) \bigr).
\]
\end{Lem}

\begin{pf}
Since $\hat\gamma_n(x)\toP\gamma(x)$, a first order Taylor
expansion yields
\[
K_{\hat\gamma_n(x)}(r_n)=K_{\gamma(x)}(r_n) +
\Lambda_n \xi_n^{(\gamma
)}(x) K'_{\tilde\gamma_n(x) }(r_n),
\]
where $\tilde\gamma_n(x)= \gamma(x)+ \Theta_n \Lambda_n \xi_n^{(\gamma
)}(x) $ with $\Theta_n\in(0,1)$.
As a consequence
\begin{eqnarray*}
\Lambda_n^{-1} \biggl( \frac{K_{\hat\gamma_n(x)}(r_n) - K_{\gamma
(x)}(r_n) }{K'_{\gamma(x)}(r_n)} \biggr)&=&
\xi_n^{(\gamma)}(x) \frac{K'_{\tilde\gamma_n(x)}(r_n)}{K'_{\gamma
(x)}(r_n)}
\\
&=& \xi_n^{(\gamma)}(x) \biggl( 1 + \frac{\int_1^{r_n} (s^{\tilde
\gamma
_n(x)-\gamma(x)}-1)s^{\gamma(x)-1}\log(s)\,\mathrm{d}s}{\int_1^{r_n} s^{\gamma
(x)-1}\log(s)\,\mathrm{d}s} \biggr).
\end{eqnarray*}
Suppose for instance $r_n\geq1$.
The assumptions yield $(\log r_n)(\tilde\gamma_n(x)-\gamma(x))\toP0$
and thus, for~$n$ large enough, with high probability,
\[
\sup_{s\in[1,r_n]}\bigl| \bigl(s^{\tilde\gamma_n(x)-\gamma(x)}-1 \bigr)\bigr|\leq2 (\log
r_n)\bigl|\tilde\gamma_n(x)-\gamma(x)\bigr| =\mathrm{o}_{\PP}(1).
\]
As a conclusion,
\[
\Lambda_n^{-1} \biggl( \frac{K_{\hat\gamma_n(x)}(r_n) - K_{\gamma
(x)}(r_n) }{K'_{\gamma(x)}(r_n)} \biggr) =
\xi^{(\gamma
)}(x) \bigl(1+\mathrm{o}_\PP(1) \bigr)
\]
and the result is proved. The case $r_n\leq1$ is easily deduced since
$K_{\gamma(x)}(1/r_n)= - K_{-\gamma(x)}(r_n)$
and $K'_{\gamma(x)}(1/r_n)= K'_{-\gamma(x)}(r_n)$.
\end{pf}

%s6.2 #&#
\subsection{Proofs of main results}
%%%%--------------------------------------------------------------------------

%pa6.2.subsubsection.1 #&#
\begin{pf*}{Proof of Theorem \protect\ref{thquant}}
Let us introduce $v_n = (nh^p\alpha_n^{-1})^{1/2}$, $\sigma_{n}(x) =
(v_nf(q(\alpha_{n}|x)|x))^{-1}$
and, for all $j=1,\ldots,J$,
\begin{eqnarray*}
W_{n,j}(x)&=&v_{n} \bigl(\habaF \bigl(q(\tau_j
\alpha_{n}|x) + \sigma_{n}(x) z_j|x \bigr)-
\baF \bigl(q(\tau_j \alpha_{n}|x) + \sigma_{n}(x)
z_j|x \bigr) \bigr),
\\
a_{n,j}(x)&=&v_{n} \bigl(\tau_j
\alpha_{n}-\baF \bigl(q(\alpha_{n,j}|x) +
\sigma_{n}(x) z_j |x \bigr) \bigr)
\end{eqnarray*}
and $z_j\in\R$. We examine the asymptotic behavior of $J$-variate function
defined by
\begin{eqnarray*}
\Phi_n(z_1,\ldots,z_J)&=&\PP \Biggl(\bigcap
_{j=1}^J \bigl\{ \sigma_{n}^{-1}(x)
\bigl(\hat q_n(\tau_j \alpha_{n}|x)-q(
\tau_j \alpha_{n}|x) \bigr)\leq z_j \bigr\}
\Biggr) \\
&=& \PP \Biggl(\bigcap_{j=1}^J \bigl
\{ W_{n,j}(x) \leq a_{n,j}(x) \bigr\} \Biggr).
\end{eqnarray*}
Let us first focus on the nonrandom term $a_{n,j}(x)$.
For each $j\in\{1,\ldots,J\}$ there exists $\theta_{n,j}\in(0,1)$
such that
\[
a_{n,j}(x) = v_{n} \sigma_{n}(x)
z_j f \bigl(q_{n,j}(x)|x \bigr) = z_j
\frac
{f(q_{n,j}(x)|x)}{f(q(\alpha_{n}|x)|x)},
\]
where
\begin{eqnarray*}
q_{n,j}(x) &=& q(\tau_j \alpha_n|x) +
\theta_{n,j}\sigma_{n}(x) z_j
\\
& = & q(\tau_j\alpha_{n}|x) + \theta_{n,j}
\frac{z_j}{\tau_j} \bigl(n h^p \alpha_n
\bigr)^{-1/2} \frac{\tau_j \alpha_n}{ f(q(\tau_j \alpha_n|x)|x)} \frac{f(q(\tau_j \alpha_n|x)|x)}{ f(q(\alpha_n|x)|x)}
\\
&=& q(\tau_j\alpha_{n}|x) + \theta_{n,j}
z_j \tau_j^{\gamma(x)} \bigl(n h^p
\alpha_n \bigr)^{-1/2} \frac{\tau_j \alpha_n}{ f(q(\tau_j \alpha_n|x)|x)} \bigl(1+
\mathrm{o}(1) \bigr),
\end{eqnarray*}
since $y\mapsto f(q(y|x)|x)$ is regularly varying at 0 with index
$\gamma(x)+1$,
see~\cite{deHaanFer}, Corollary~1.1.10, equation (1.1.33). Now, in view
of~\cite{deHaanFer}, Theorem~1.2.6
and~\cite{deHaanFer}, Remark~1.2.7, a possible choice of the auxiliary
function is
%
%e16 #&#
\renewcommand{\theequation}{\arabic{equation}}
\begin{equation}
\label{fctaux} a(t|x)=\frac{\baF(t|x)}{f(t|x)} \bigl(1+\mathrm {o}(1) \bigr),
\end{equation}
leading to
\[
q_{n,j}(x)=q(\tau_j\alpha_{n}|x) +
\theta_{n,j} z_j \tau_j^{\gamma(x)} \bigl(n
h^p \alpha_n \bigr)^{-1/2} a \bigl(q(
\tau_j \alpha_n|x)|x \bigr) \bigl(1+\mathrm{o}(1) \bigr).
\]
Applying Lemma~\ref{lemunif} with $z_n(x)=q(\tau_j\alpha_{n}|x)$,
$t_n(x)=\theta_{n,j} z_j \tau_j^{\gamma(x)} (n h^p \alpha_n)^{-1/2}(1+\mathrm{o}(1))$
and $t_0(x)=0$ yields
\[
\frac{\baF(q_{n,j}(x)|x)}{\tau_j\alpha_n} \to K_{\gamma
(x)}^{-1}(0) = 1
\]
as $n\to\infty$. Recalling that $y\mapsto f(q(y|x)|x)$ is regularly varying,
we have
\[
\frac{f(q_{n,j}(x)|x)}{f(q(\alpha_{n}|x)|x)} \to\tau_j^{\gamma(x)+1}
\]
as $n\to\infty$ and therefore
%
%e17 #&#
\begin{equation}
\label{anj} a_{n,j}(x) = z_j \tau_j^{\gamma(x)+1}
\bigl(1+\mathrm{o}(1) \bigr), \qquad j=1,\ldots,J.
\end{equation}
Let us now turn to the random term $W_{n,j}(x)$. Let us define
$z_{n,j}(x)=q(\tau_j\alpha_{n}|x)+\sigma_{n}(x)z_j$ for $j=1,\ldots,J$,
$y_n(x)=q(\alpha_n|x)$, and consider the expansion
\[
\frac{z_{n,j}(x)-y_n(x)}{a(y_n(x)|x)} = \frac{q(\tau_j\alpha_n|x)-q(\alpha_n|x)}{a(q(\alpha_n|x)|x)} + \frac{\sigma_{n}(x)z_j}{a(q(\alpha_n|x)|x)} .
\]
From~(\ref{F1quant}), we have
\[
\lim_{n \to\infty} \frac{q(\tau_j\alpha_n|x)-q(\alpha_n|x)}{a(q(\alpha_n|x)|x)} = K_{\gamma(x)}(1/
\tau_j),
\]
and
\[
\lim_{n \to\infty} \frac{\sigma_{n}(x)z_j}{a(q(\alpha_n|x)|x)} = 0,
\]
leading to $z_{n,j}(x)= y_n(x) + K_{\gamma(x)}(1/\tau_j)a(y_n(x)|x)(1+\mathrm{o}(1))$. Introducing $\beta_{n,j}(x)= \baF
(y_{n,j}(x)|x)$,
the oscillation $\omega_n(x)$ can be rewritten as
\[
\omega_n(x)=\max_{j=1,\ldots,J}\sup_{x'\in B(x,h)} \biggl\llvert
\frac
{\baF
(q(\beta_{n,j}(x)|x)|x')}{\beta_{n,j}(x)}-1 \biggr\rrvert.
\]
For all $\kappa\in(0,\tau_J)$ and $j=1,\ldots,J$, we eventually have
$z_{n,j}(x)\in[y_n(x),z_n(x)]$ where $z_n(x):= y_n(x) + K_{\gamma
(x)}(2/\kappa)a(y_n(x)|x)$
and thus $\beta_{n,j}(x)\in[\baF(z_n(x)|x), \alpha_n]$ eventually.
Now, Lemma~\ref{lemunif} implies that $\baF(z_n(x)|x)/\alpha_n \to
\kappa/2$ as $n\to\infty$ and thus, for $n$ large enough,
$\beta_{n,j}(x)\in[\kappa\alpha_n, \alpha_n]$. Consequently,
$\omega_n(x)\leq\Delta_{\kappa}(\alpha_n,x)$.
Applying Proposition~\ref{thproba} and Lemma~\ref{lemunif} yields
\[
W_{n,j}(x) = \frac{\baF(z_{n,j}(x)|x)}{\alpha_n} \xi_{n,j} ={
\tau_j}\xi_{n,j} \bigl(1+\mathrm{o}(1) \bigr),
\]
where $\xi_n=(\xi_{n,1},\ldots,\xi_{n,j})^t$ converges to a centered
Gaussian random vector with covariance matrix $\|K\|_2^2/g(x) V$.
Taking into account of~(\ref{anj}), the results follows.
\end{pf*}

%pa6.2.subsubsection.2 #&#
\begin{pf*}{Proof of Corollary \protect\ref{coroq}} Let us remark that,
from~(\ref{fctaux}),
\begin{eqnarray*}
&& \Bigl\{ f \bigl(q(\alpha_n|x)|x \bigr) \sqrt{nh^p
\alpha_n^{-1}} \bigl(\hat q_n(
\tau_j \alpha_n|x)-q(\tau_j
\alpha_n|x) \bigr) \Bigr\}_{j=1,\ldots,J}
\\
&&\quad= \biggl\{ q(\alpha_n|x)\frac{f(q(\alpha_n|x)|x)}{\alpha_n} \frac
{q(\tau_j\alpha_n|x)}{q(\alpha_n|x)}
\sqrt{nh^p\alpha_n} \biggl(\frac{\hat
q_n(\tau_j \alpha_n|x)}{q(\tau_j \alpha_n|x)}-1 \biggr)
\biggr\}_{j=1,\ldots,J}
\\
&&\quad= \biggl\{ \frac{q(\alpha_n|x)}{a(q(\alpha_n|x)|x)} \frac{q(\tau_j\alpha_n|x)}{q(\alpha_n|x)} \sqrt{nh^p
\alpha_n} \biggl(\frac{\hat
q_n(\tau_j \alpha_n|x)}{q(\tau_j \alpha_n|x)}-1 \biggr) \biggr\}_{j=1,\ldots,J}
\bigl(1+\mathrm{o}(1) \bigr)
\\
&&\quad= \biggl\{ \frac{q(\alpha_n|x)}{a(q(\alpha_n|x)|x)} \tau_j^{-(\gamma
(x)\vee0)}
\sqrt{nh^p\alpha_n} \biggl(\frac{\hat q_n(\tau_j \alpha_n|x)}{q(\tau_j \alpha_n|x)}-1 \biggr)
\biggr\}_{j=1,\ldots,J} \bigl(1+\mathrm{o}(1) \bigr),
\end{eqnarray*}
in view of Lemma~\ref{lemrapquant}. The result follows from
Theorem~\ref
{thquant}.
\end{pf*}

%pa6.2.subsubsection.3 #&#
\begin{pf*}{Proof of Theorem \protect\ref{thquant2}} By definition,
\[
q_n(\beta_n|x)=q_n(\alpha_n|x) +
\bigl(K_{\gamma(x)}(\alpha_n/\beta_n) + b(
\beta_n/\alpha_n,\alpha_n) \bigr) a \bigl(q(
\alpha_n|x)|x \bigr)
\]
and thus,
the following expansion can be easily established:
\begin{eqnarray*}
&&\Lambda_n^{-1} \biggl(\frac{\tilde q_n(\beta_n|x)-q(\beta_n|x)}{a(q(\alpha_n|x)|x) K'_{\gamma(x)}(\alpha_n/\beta_n)} \biggr)\\
&&\quad=
\Lambda_n^{-1} \biggl(\frac{\hat q_n(\alpha_n|x)-q(\alpha_n|x)}{a(q(\alpha_n|x)|x) K'_{\gamma(x)}(\alpha_n/\beta_n)} \biggr)
\\
&&\qquad{}+ \Lambda_n^{-1} \biggl(\frac{K_{\hat\gamma_n(x)}(\alpha_n/\beta_n)-K_{\gamma(x)}(\alpha_n/\beta_n)} {K'_{\gamma
(x)}(\alpha_n/\beta_n)} \biggr)
\frac{\hat a_n(x)}{a(q(\alpha_n|x)|x)}
\\
&&\qquad{}+ \Lambda_n^{-1} \frac{K_{\gamma(x)}(\alpha_n/\beta_n)} {
K'_{\gamma
(x)}(\alpha_n/\beta_n)} \biggl(
\frac{\hat a_n(x)}{a(q(\alpha_n|x)|x)}-1 \biggr)
\\
&&\qquad{}+ \Lambda_n^{-1} \frac{b(\beta_n/\alpha_n,\alpha_n)}{K'_{\gamma
(x)}(\alpha_n/\beta_n)}
\\
&&\quad=: T_{n,1} + T_{n,2} + T_{n,3} +
T_{n,4}.
\end{eqnarray*}
Introducing
\[
\bigl(\xi_n^{(\gamma)}(x), \xi_n^{(a)}(x),
\xi_n^{(q)}(x) \bigr):= \Lambda_n^{-1}
\biggl( \hat\gamma_n(x)-\gamma(x), \frac{\hat a_n(x)}{a(q(\alpha_n|x)|x)}-1,
\frac{\hat q_n( \alpha_n|x)-q( \alpha_n|x)}{a(q(\alpha_n|x)|x)} \biggr),
\]
from and remarking that, when $u\to\infty$,
%
%e18 #&#
\renewcommand{\theequation}{\arabic{equation}}
\begin{equation}
\label{limites} K'_z(u)= \bigl(1+\mathrm{o}(1) \bigr)
\cases{\displaystyle\frac{1}{z^2}, &\quad $\mbox{if } z<0,$
\vspace*{2pt}\cr
\displaystyle\frac{\log^2(u)}{2}, & \quad$\mbox{if } z=0,$
\vspace*{2pt}\cr
\displaystyle\frac{u^z\log(u)}{z}, &\quad $\mbox{if } z>0,$}
\end{equation}
the first term can be rewritten as
%
%e19 #&#
\begin{equation}
\label{tn1} T_{n,1} = \frac{\xi_n^{(q)}(x)}{K'_{\gamma(x)}(\alpha_n/\beta_n)} = \bigl(\gamma(x)\wedge0
\bigr)^2 \xi_n^{(q)}(x) \bigl(1 +
\mathrm{o}_{\PP}(1) \bigr).
\end{equation}
Second, $\Lambda_n\to0$ and (\ref{condtriplet}) entail ${\hat
a_n(x)}/{a(q(\alpha_n|x)|x)}\toP1$ and thus
%
%e20 #&#
\begin{equation}
\label{tn2} T_{n,2} = \Lambda_n^{-1} \biggl(
\frac{K_{\hat\gamma_n(x)}(\alpha_n/\beta_n)-K_{\gamma
(x)}(\alpha_n/\beta_n)} {K'_{\gamma(x)}(\alpha_n/\beta_n)} \biggr) \bigl(1+\mathrm{o}_\PP(1) \bigr) =
\xi_n^{(\gamma)}(x) \bigl(1+\mathrm{o}_\PP(1) \bigr),
\end{equation}
from Lemma~\ref{lemNAK}.
From (\ref{condtriplet}), (\ref{limites}), and in view of
\[
K_z(u)= \bigl(1+\mathrm{o}(1) \bigr) %
\cases{\displaystyle-\frac{1}{z}, &\quad $\mbox{if } z<0,$
\vspace*{2pt}\cr
\log(u), & \quad $\mbox{if } z=0,$
\vspace*{2pt}\cr
\displaystyle\frac{u^z}{z}, &\quad  $\mbox{if } z>0.$ }
\]
The third term can be rewritten as
%
%e21 #&#
\begin{equation}
\label{tn3} T_{n,3} = \xi_n^{(a)}(x)
\frac{K_{\gamma(x)}(\alpha_n/\beta_n)}{K'_{\gamma(x)}(\alpha_n/\beta_n)} = - \bigl(\gamma(x)\wedge0 \bigr) \xi_n^{(a)}(x)
\bigl(1 + \mathrm{o}_{\PP}(1) \bigr).
\end{equation}
Finally, $T_{n,4}=\mathrm{o}_{\PP}(1)$ by assumption and the
conclusion follows
from (\ref{condtriplet}), (\ref{tn1}), (\ref{tn2}) and~(\ref
{tn3}).
\end{pf*}

%pa6.2.subsubsection.4 #&#
\begin{pf*}{Proof of Theorem \protect\ref{triplet}} The proof consists
in deriving asymptotic expansions
for the three considered random variables.
(i) Let us first introduce
%
%e22 #&#
\renewcommand{\theequation}{\arabic{equation}}
\begin{equation}
\label{gammanj} \gamma_{n,j}(x) = \frac{1}{\log r} \log \biggl(
\frac{\hat q_n(\tau_j \alpha_n|x) - \hat q_n(\tau_{j+1}\alpha_n|x)} {
\hat q_n(\tau_{j+1} \alpha_n|x) - \hat q_n(\tau_{j+2}\alpha_n|x)} \biggr)
\end{equation}
such that $\hat{\gamma}_n^{\mathrm{RP}}(x)= \sum_{j=1}^{J-2} \pi_j
\gamma_{n,j}(x)$.
From Theorem~\ref{thquant} and in view of~(\ref{F1quant}), we have, for
all $j=1,\ldots,J$,
\[
\hat q_n(\tau_j \alpha_n|x) = q(
\alpha_n|x) + a \bigl(q(\alpha_n|x)|x \bigr) \bigl(
K_{\gamma(x)}(1/\tau_j) + b(\tau_j,
\alpha_n|x) \bigr) + \sigma_n(x) \xi_{j,n},
\]
with $\sigma_n^{-1}(x)=f(q(\alpha_n|x)|x) \sqrt{nh^p\alpha_n^{-1}}$ and
where the random vector $\xi_n=(\xi_{j,n})_{j=1,\ldots,J}$ is asymptotically
Gaussian, centered, with covariance matrix
$
\|K\|_2^2/g(x) \Sigma(x)
$.
Introducing
\begin{eqnarray*}
\eta_n(x)&:=& \max_{j=1,\ldots,J} \bigl| b(\tau_j,
\alpha_n|x) \bigr|,
\\
\varepsilon_n&:=&\sigma_n(x) / a \bigl(q(
\alpha_n|x)|x \bigr)= \bigl(n h^p\alpha_n
\bigr)^{-1/2} \bigl(1+\mathrm{o}(1) \bigr),
\end{eqnarray*}
see~\cite{deHaanFer}, it follows that
%
%e23 #&#
\begin{eqnarray}\label{Devquant2}
&&
\frac{\hat q_n(\tau_j \alpha_n|x)-\hat q_n(\tau_{j+1} \alpha_n|x)}{a(q(\alpha_n|x)|x)}\nonumber\\
&&\quad =\varepsilon_n(\xi_{j,n}-
\xi_{j+1,n})
\nonumber
\\[-8pt]
\\[-8pt]
\nonumber
&&\qquad{}+ K_{\gamma(x)}(1/\tau_j)-K_{\gamma(x)}(1/
\tau_{j+1}) + b(\tau_j,\alpha_n|x)- b(
\tau_{j+1},\alpha_n|x)
\\
 &&\quad= \varepsilon_n(\xi_{j,n}-
\xi_{j+1,n}) + K_{\gamma(x)}(r) r^{-\gamma
(x)j} + \mathrm{O} \bigl(
\eta_n(x) \bigr).\nonumber
\end{eqnarray}
Replacing in (\ref{gammanj}), we obtain
\[
(\log r) \gamma_{n,j}(x) = \log \biggl( \frac{\varepsilon_n(\xi_{j,n}-\xi_{j+1,n}) + K_{\gamma(x)}(r)
r^{-\gamma
(x)j} + \mathrm{O}( \eta_n(x))} {
\varepsilon_n(\xi_{j+1,n}-\xi_{j+2,n}) + K_{\gamma(x)}(r)
r^{-\gamma
(x)(j+1)} + \mathrm{O}( \eta_n(x))} \biggr),
\]
or equivalently,
\begin{eqnarray*}
(\log r) \bigl(\gamma_{n,j}(x)-\gamma(x) \bigr) &= & \log \biggl( 1+
\frac{\varepsilon_n(\xi_{j,n}-\xi_{j+1,n})r^{\gamma
(x)j}}{ K_{\gamma(x)}(r)}+ \mathrm{O} \bigl( \eta_n(x) \bigr) \biggr)
\\
&&{}-\log \biggl( 1+ \frac{\varepsilon_n(\xi_{j+1,n}-\xi_{j+2,n})r^{\gamma
(x)(j+1)}}{ K_{\gamma(x)}(r)}+ \mathrm{O} \bigl( \eta_n(x) \bigr)
\biggr).
\end{eqnarray*}
A first order Taylor expansion yields
\begin{eqnarray*}
&&(\log r) \varepsilon_n^{-1} \bigl(\gamma_{n,j}(x)
-\gamma(x) \bigr)\\
&&\quad = \frac{r^{\gamma(x)j}}{K_{\gamma(x)}(r)} \bigl( \xi_{j,n} -
\bigl(1+r^{\gamma(x)} \bigr) \xi_{j+1,n} + r^{\gamma(x)}
\xi_{j+2,n} \bigr)+\mathrm{O}\bigl(\varepsilon_n^{-1}
\eta_n(x)\bigr) + \mathrm{o}_{\PP}(1)
\end{eqnarray*}
and thus, under the assumption $(nh^p\alpha_n)^{1/2} \eta_n(x)\to0$ as
$n\to\infty$,
\begin{eqnarray*}
&&\sqrt{nh^p\alpha_n} \bigl(\hat{\gamma}_n^{\mathrm{RP}}(x)-
\gamma (x) \bigr) \\
&&\quad= \frac{1}{(\log r) K_{\gamma(x)}(r)} \sum_{j=1}^{J-2}
\pi_j r^{\gamma(x)j} \bigl( \xi_{j,n} -
\bigl(1+r^{\gamma
(x)} \bigr) \xi_{j+1,n} + r^{\gamma(x)}
\xi_{j+2,n} \bigr)+ \mathrm{o}_{\PP}(1).
\end{eqnarray*}
Defining for the sake of simplicity
$\pi_{-1}=\pi_0=\pi_{J-1}=\pi_{J}=0$, $\beta_0^{(\gamma)}=\frac
{1}{\log
r}$, $\beta_1^{(\gamma)}=-\frac{1+r^{-\gamma(x)}}{\log(r)}$ and
$\beta_2^{(\gamma)}=\frac{r^{-\gamma(x)}}{\log(r)}$,
we end up with
%
%e24 #&#
\begin{eqnarray}\label{devgamma}
\xi_n^{(\gamma)}(x)&:= & \sqrt{nh^p
\alpha_n} \bigl(\hat{\gamma}_n^{\mathrm{RP}}(x)-\gamma(x)
\bigr)
 \nonumber
 \\[-8pt]
 \\[-8pt]
 \nonumber
 &=& \frac{1}{K_{\gamma(x)}(r)} \sum_{j=1}^J
r^{\gamma(x)j} \bigl( \beta_0^{(\gamma)} \pi_j +
\beta_1^{(\gamma)} \pi_{j-1} + \beta_2^{(\gamma
)}
\pi_{j-2} \bigr)\xi_{j,n} + \mathrm{o}_{\PP}(1).
\end{eqnarray}
(ii) Second, let us now consider
\[
a_{n,j}(x) = \frac{r^{\hat{\gamma}_n^{\mathrm{RP}}(x)j}(\hat
q_n(\tau_j \alpha_n|x) - \hat q_n(\tau_{j+1}\alpha_n|x))} {
K_{\hat{\gamma}_n^{\mathrm{RP}}(x)}(r)}
\]
such that $\hat a_n(x) = \sum_{j=1}^{J-2} \pi_j a_{n,j}(x)$.
From~(\ref{Devquant2}), it follows that, for all $j=1,\ldots,J$,
\[
\frac{a_{n,j}(x)}{a(q(\alpha_n|x)|x)} = \frac{r^{\hat{\gamma
}_n^{\mathrm{RP}}(x)j}}{
K_{\hat{\gamma}_n^{\mathrm{RP}}(x)}(r)} \bigl( \varepsilon_n (
\xi_{j,n}-\xi_{j+1,n}) + K_{\gamma(x)}(r) r^{-\gamma(x)j} +
\mathrm{O} \bigl(\eta_n(x) \bigr) \bigr).
\]
Remarking that $\hat{\gamma}_n^{\mathrm{RP}}(x)=\gamma(x) +
(nh^p\alpha_n)^{-1/2}\xi_n^{(\gamma)}(x)$, Lemma~\ref{lemNAK} yields
\begin{eqnarray*}
\frac{a_{n,j}(x)}{a(q(\alpha_n|x)|x)} &=& \frac{ 1+ ({r^{\gamma
(x) j}}/{K_{\gamma(x)}(r)}) \varepsilon_n
(\xi_{j,n}-\xi_{j+1,n}) + \mathrm{O}(\eta_n(x))} {
1+ ({K'_{\gamma(x)}(r)}/{K_{\gamma(x)}(r)}) (nh^p\alpha_n)^{-1/2}
\xi_n^{(\gamma)}(x)(1+\mathrm{o}_{\PP}(1))}\\
&&{}\times  \exp \bigl( \xi_n^{(\gamma)}(x)
j\log(r) \bigl(nh^p\alpha_n \bigr)^{-1/2} \bigr).
\end{eqnarray*}
A first order Taylor expansion thus gives
\begin{eqnarray*}
\sqrt{nh^p\alpha_n} \biggl(\frac{a_{n,j}(x)}{a(q(\alpha_n|x)|x)}-1 \biggr)
&= &\xi_n^{(\gamma)}(x) \biggl(j\log(r)- \frac{K'_{\gamma
(x)}(r)}{K_{\gamma(x)}(r)}
\biggr) + \frac{r^{\gamma(x) j}}{K_{\gamma(x)}(r)} (\xi_{j,n}-\xi_{j+1,n})
\\
&&{} + \mathrm{O} \bigl(\sqrt{nh^p\alpha_n}\eta_n(x)
\bigr) + \mathrm{o}_{\PP}(1).
\end{eqnarray*}
Recalling that $\pi_{-1}=\pi_0=\pi_{J-1}=\pi_{J}=0$ and introducing
\begin{eqnarray*}\E(\pi)&=&\sum_{j=1}^{J}
j\pi_j, \qquad \beta_1^{(a)} = -r^{-\gamma(x)} -
\bigl(r^{-\gamma(x)}+1 \bigr) \biggl( \E(\pi) - \frac{K'_{\gamma
(x)}(r)}{\log(r)
K_{\gamma(x)}(r)} \biggr),
\\
\beta_0^{(a)}&=&1+\E(\pi) - \frac{K'_{\gamma(x)}(r)}{\log(r)
K_{\gamma(x)}(r)}, \qquad
\beta_2^{(a)} = r^{-\gamma(x)} \biggl( \E(\pi) -
\frac
{K'_{\gamma
(x)}(r)}{\log(r) K_{\gamma(x)}(r)} \biggr)
\end{eqnarray*}
it follows that
%
%e25 #&#
\begin{eqnarray}\label{deva}
\nonumber
\xi_{n}^{(a)}(x)&:= & \sqrt{nh^p
\alpha_n} \biggl( \frac
{a_n(x)}{a(q(\alpha_n|x)|x)} -1 \biggr)
\\
&= & \biggl( \E(\pi) \log(r) - \frac{K'_{\gamma(x)}(r)}{K_{\gamma
(x)}(r)} \biggr)
\xi_n^{(\gamma)(x)}
\nonumber
\\[-8pt]
\\[-8pt]
\nonumber
&&{} + \frac{1}{K_{\gamma(x)}(r)} \sum
_{j=1}^{J} r^{\gamma(x)j} \bigl(\pi_j
-r^{-\gamma(x)} \pi_{j-1} \bigr)\xi_{j,n} +
\mathrm{o}_{\PP}(1)
\\
 &=& \frac{1}{K_{\gamma(x)}(r)} \sum_{j=1}^{J}
r^{\gamma(x)j} \bigl(\beta_0^{(a)} \pi_j+
\beta_1^{(a)}\pi_{j-1}+\beta_2^{(a)}
\pi_{j-2} \bigr) \xi_{j,n} + \mathrm{o}_{\PP}(1)\nonumber
\end{eqnarray}
in view of (\ref{devgamma}).
(iii) Third, Corollary~\ref{coroq} states that
%
%e26 #&#
\begin{equation}
\label{devq} \xi_{1,n}= \frac{\sqrt{nh^p\alpha_n}}{a(q(\alpha_n|x)|x)} \bigl(\hat
q_n( \alpha_n|x)- q( \alpha_n|x) \bigr)
\end{equation}
is asymptotically Gaussian.
Finally, collecting (\ref{devgamma}), (\ref{deva}) and (\ref{devq}),
\[
\bigl(\xi_n^{(\gamma)}(x),\xi_n^{(a)}(x),
\xi_{1,n} \bigr)^t = \frac
{1}{K_{\gamma
(x)}(r)} A(x)
\xi_n + \mathrm{o}_{\PP}(1),
\]
where $A(x)$ is the $3\times J$ matrix defined by
\begin{eqnarray*}
A_{1,j}(x) & = &r^{\gamma(x)j} \bigl( \beta_0^{(\gamma)}
\pi_j + \beta_1^{(\gamma)} \pi_{j-1} +
\beta_2^{(\gamma)}\pi_{j-2} \bigr)
\\
A_{2,j}(x) & = &r^{\gamma(x)j} \bigl(\beta_0^{(a)}
\pi_j+\beta_1^{(a)}\pi_{j-1}+
\beta_2^{(a)}\pi_{j-2} \bigr)
\\
A_{3,j}(x) &=& K_{\gamma(x)}(r) \I\{j=1\}
\end{eqnarray*}
for all $j=1,\ldots,J$.
It is thus clear
that the random vector $(\xi_n^{(\gamma)}(x),\xi_n^{(a)}(x),\xi_{1,n})^t $ converges in distribution to
a centered Gaussian random vector with covariance matrix
%
%e27 #&#
\begin{equation}
\label{matS} \frac{\|K\|_2^2}{g(x)K^2_{\gamma(x)}(r)} A(x) \Sigma (x) A^t(x)=:S(x).
\end{equation}
\end{pf*}
\end{appendix}

\section*{Acknowledgements}
The authors thank the Editor
and both anonymous Associate Editor and reviewer for their valuable
feedback on this paper.
%%%%%%%%%%%%%%%%%%%%%%%%%

% imsref loaded by akundreckaite, 2012-10-02 11:19:22
%

%suskaldyti doi

\printhistory

\end{document}